\documentclass{article} 
\usepackage{epsfig,epic,eepic,latexsym, amssymb}
\input xy

\xyoption{all}
\xyoption{all}
\makeindex

\def \pf{PROOF: }
\def \QED{\hfill\hbox{\hskip 4pt
                \vrule width 5pt height 6pt depth 1.5pt}}
\def \epf{\QED\\}
\def \spa{$\;\;\;\;$}

\newcommand{\Cl}{ \{   \hskip -3pt \mid  }
\newcommand{\Cr}{ \mid \hskip -3pt \}  }
\newcommand{\Sl}{ [ \hskip -1.5pt [  }
\newcommand{\Sr}{ ] \hskip -1.5pt ]  }

\newcommand{\V}{  {\cal V} }
\newcommand{\A}{  {\cal A} }
\newcommand{\B}{  {\cal B} }
\newcommand{\C}{  {\cal C} }
\newcommand{\D}{  {\cal D} }
\newcommand{\K}{  {\cal K} }
\newcommand{\LL}{  {\cal L} }
\newcommand{\N}{ {\cal N} }
\newcommand{\Q}{ {\cal Q} }

\newcommand{\VCat}{ {\cal V}\mbox{-}\mathit{\bf Cat}}
\newcommand{\VCAT}{ {\cal V}\mbox{-}\mathit{\bf CAT}}
\newcommand{\VMod}{ {\cal V}\mbox{-}\mathit{\bf Mod}}

\newcommand{\Set}{ \mathit{\bf Set}}
\newcommand{\Cat}{ \mathit{\bf Cat}}
\newcommand{\CAT}{ \mathit{\bf CAT}}

\newcommand{\SCat}{ \Set\mbox{-}\mathit{\bf Cat} }
\newcommand{\SCAT}{ \Set\mbox{-}\mathit{\bf CAT} }

\newcommand{\I}{{ \cal I} }

\newcommand{\Co}{\mathit{\bf Cocts}}
\newcommand{\Cont}{\mathit{\bf Conts}}

\newcommand{\DCont}{ \D\mbox{-}\mathit{\bf Conts}}

\newcommand{\PhiCo}{\Phi\mbox{-}\mathit{\bf Cocts}}
\newcommand{\PhiCont}{\Phi\mbox{-}\mathit{\bf Conts}}
\newcommand{\PhiSCont}{{\Phi}^*\mbox{-}\mathit{\bf Conts}}

\newcommand{\QCo}{\Q\mbox{-}\mathit{\bf Cocts}}
\newcommand{\PsiCo}{\Psi\mbox{-}\mathit{\bf Cocts}}
\newcommand{\PsiCont}{\Psi\mbox{-}\mathit{\bf Conts}}

\newcommand{\ZCont}{ 0\mbox{-}\mathit{\bf Conts}}

\newcommand{\p}{{\cal P}}
\newcommand{\PCo}{{\cal P}\mbox{-}\mathit{\bf Cocts}}
\newcommand{\PCont}{{\cal P}\mbox{-}\mathit{\bf Conts}}

\newtheorem{theorem}{Theorem}[section]
\newtheorem{definition}[theorem]{Definition}
\newtheorem{proposition}[theorem]{Proposition}
\newtheorem{lemma}[theorem]{Lemma}
\newtheorem{example}[theorem]{Example}
\newtheorem{corollary}[theorem]{Corollary}
\newtheorem{remarks}[theorem]{Remarks}
\newtheorem{remark}[theorem]{Remark}
\newtheorem{fact}[theorem]{} 
\newtheorem{observation}[theorem]{Observation}

\title{Notes on enriched categories with colimits of some class
 }
\author{G.M. Kelly\thanks{
School of Mathematics and Statistics F07,
University of Sydney, NSW 2006
Australia, maxk@maths.usyd.edu.au} \thanks{The first author gratefully
acknowledges the support of the Australian Research Council.}
and V.Schmitt\thanks{The University of Leicester,
School of Mathematics and Computer Science,
LE17RH Leicester United Kingdom, vs27@mcs.le.ac.uk}
\thanks{Grant GR/R63004/01 of the Engineering and Physical Sciences 
Research Council gratefully acknowledged.} }

\begin{document}
\maketitle                      

\begin{abstract}
The paper is in essence a survey
of categories having $\phi$-weighted colimits 
for all the weights $\phi$ in some class $\Phi$.
We introduce the class $\Phi^+$ of {\em $\Phi$-flat} 
weights which are those $\psi$ for which $\psi$-colimits commute 
in the base $\V$ with limits having weights in $\Phi$; and 
the class $\Phi^-$ of {\em $\Phi$-atomic}
weights, which are those $\psi$ for which $\psi$-limits commute 
in the base $\V$ with colimits having weights in $\Phi$.
We show that both these classes are {\em saturated} 
(that is, what was called
{\em closed} in the terminology of \cite{AK88}).
We prove that for the class $\p$ of {\em all} weights,
the classes 
$\p^+$ and $\p^-$ both coincide with the class
$\Q$ of {\em absolute} weights.
For any class $\Phi$ and any category $\A$, we 
have the free $\Phi$-cocompletion $\Phi(\A)$ of $\A$;
and we recognize $\Q(\A)$ as the Cauchy-completion 
of $\A$. We study the equivalence between ${(\Q(\A^{op}))}^{op}$
and $\Q(\A)$, which we exhibit as the restriction 
of the Isbell adjunction between ${[\A,\V]}^{op}$
and $[\A^{op},\V]$ when $\A$ is small; and we give a new Morita theorem
for any class $\Phi$ containing $\Q$.
We end with the study of $\Phi$-continuous weights
and their relation to the $\Phi$-flat weights.
\end{abstract}

\begin{section}{Introduction}\label{section1}
The present observations had their beginnings 
in an analysis of the results obtained by Borceux, 
Quintero and Rosick\'y in their article
\cite{BQR98}, which in turn followed on from that of 
Borceux and Quintero \cite{BQ96}. These authors 
were concerned with extending
to the enriched case the notion of accessible category
and its properties, described for ordinary categories
in the books \cite{MP89} of Makkai and Par\'e and 
\cite{AR94} of Ad\`amek and Rosick\'y. They were led to discuss
categories -- now meaning $\V$-categories -- with finite 
limits (in a suitable sense), or more generally with 
$\alpha$-small limits, or with filtered colimits (in a suitable
sense), and more generally with $\alpha$-filtered colimits,
or again with $\alpha$-flat colimits,
and to discuss the connexions between these classes
of limits and of colimits. When we looked in detail at their work,
we observed that many of the properties they discussed
hold in fact for categories having colimits of {\em any}
given class $\Phi$, while others hold when $\Phi$ is the class
of colimits commuting in the base category $\V$ with the limits
of some class $\Psi$ -- such particular properties as finiteness
or filteredness arising only as special cases of the {\em general} 
results. Approaching in this abstract way, not generalizations of 
accessible categories as such, but the study of categories with colimits 
(or limits) of some class, brings considerable notional 
simplifications.\\ 

Although our original positive results are limited 
in number, their value may be judged by the extra light they 
cast on several of the results in \cite{BQR98}.
To expound these results, it has seemed to us necessary
to repeat some known facts so as to provide the proper context.
The outcome is that we have produced a rather complete study
of categories having colimits of a given class, which is
to a large extent self-contained: a kind of survey paper
containing a fair number of original results.\\

We begin by reviewing and completing some known material in 
the first sections: in Section \ref{section2} the general notions 
of weighted limits and colimits for enriched categories; in Section 
\ref{section3} the free $\Phi$-cocompletion $\Phi(\A)$ of a
$\V$-category $\A$; and in Section \ref{section4} results on the 
recognition of categories of the form $\Phi(\A)$.\\

Section \ref{section5} treats generally the commutation of limits 
and colimits in the base $\V$: it introduces classes of the form 
$\Phi^+$ of {\em $\Phi$-flat} weights -- those weights whose 
colimits in $\V$ commute with $\Phi$-weighted limits -- and 
classes of the form $\Phi^-$ of {\em $\Phi$-atomic} 
weights -- those weights whose limits in $\V$ commute 
with $\Phi$-weighted colimits. We show that each of these
classes is saturated.\\
 
Section \ref{section6} focuses on the class $\Q = \p^-$ where 
$\p$ is the class of {\em all} (small) weights;
this $\Q$ is the class of {\em small projective}
or {\em atomic} weights, which is also, as Street showed
in \cite{Str83}, the class of {\em absolute} weights.
We show that $\Q$ is also the class $\p^+$ of 
$\p$-flat weights.  
We recall that a weight $\phi: \K^{op} \rightarrow \V$
corresponds to a module 
$\overline{\phi}: \xymatrix{ \I \ar[r]|{\circ} & \K}$,
while a weight $\psi: \K \rightarrow \V$ corresponds
to a module $\underline{\psi} : 
\xymatrix{ \K \ar[r]|{\circ} & \I}$;
and we recall that the relation between a left adjoint
module $\overline{\phi}$ and its right adjoint 
$\underline{\psi}$ gives rise to an equivalence 
between $(\Q(\K^{op}))^{op}$ and $\Q(\K)$,
which is in fact the restriction to the small 
projectives of the {\em Isbell Adjunction}
between $[\K,\V]^{op}$ and $[\K^{op},\V]$.\\
 
Section \ref{section7} studies the Cauchy-completion $\Q(\A)$
for a general category $\A$ and gives an extension 
of the classical Morita theorem:
for any class $\Phi$ containing $\Q$ we have 
$\Phi(\A) \simeq \Phi(\B)$ if and only if
$\Q(\A) \simeq \Q(\B)$. (We use $\cong$ to denote isomorphism
and $\simeq$ to denote equivalence.)\\

Finally we consider in section \ref{section8} the class 
of $\Phi$-continuous functors $\N^{op} \rightarrow \V$,
where $\N$ is a small category admitting $\Phi$-colimits;
and we compare these with the $\Phi$-flat functors.
For $\V = \Set$, some special cases of the results 
here appeared in \cite{ABLR02}.\\

We have benefited greatly from discussions with 
Francis Borceux and with Ross Street, both of whom 
have contributed significantly to the improvement of
our exposition; 
we thankfully acknowledge their help.
\end{section}

\begin{section}{Revision of weighted limits and colimits}
\label{section2} 
The necessary background knowledge about 
enriched categories is largely contained in \cite{Kel82}, 
augmented by \cite{Kel82-2} and the Albert-Kelly article 
\cite{AK88}.\\

We deal with categories enriched in a symmetric monoidal closed 
category $\V$, supposing as usual that the ordinary category $\V_0$ 
underlying $\V$ is locally small, complete and cocomplete. 
(A set is {\em small} when its cardinal is less than a chosen inaccessible 
cardinal $\infty$, and a category is {\em locally small} when each
of its hom-sets is small.) We henceforth use ``category'', ``functor'', and 
``natural transformation'' to mean ``$\V$-category'',
``$\V$-functor'', and ``$\V$-natural
transformation'', except when more precision is needed. 
We call a $\V$-category {\em small} when its set of 
{\em isomorphism classes} of objects is a small set;
a $\V$-category that is not small is sometimes said to be {\em large}.
$\VCAT$ is the 2-category of $\V$-categories, whereas 
$\VCat$ is that of small $\V$-categories. $\Set$ is the category of small
sets, $\Cat = \SCat$ is the 2-category of small categories, and
$\CAT = \SCAT$ is the 2-category of locally small categories.\\  

A {\em weight} is a functor $\phi: \K^{op} \rightarrow \V$ 
with domain $\K^{op}$ {\em small}; weights were called {\em indexing-types}
in \cite{Kel82}, \cite{Kel82-2} and \cite{AK88}, where weighted limits 
were called {\em indexed limits}. (A functor with codomain $\V$ is often 
called a {\em presheaf}; so that a weight is a presheaf with a small domain.)
Recall that the {\em $\phi$-weighted
limit} $\{\phi, T\}$ of a functor $T:\K^{op} \rightarrow \A$ is defined
representably by
\begin{fact}\label{fac2.1}\spa
$\A(a, \{ \phi, T \}) \cong [\K^{op}, \V](\phi, \A(a,T-))$,
\end{fact}
while the {\em $\phi$-weighted colimit} $\phi*S$ of 
$S:\K \rightarrow \A$ is defined dually by
\begin{fact}\label{fac2.2}\spa
$\A(\phi * S,a) \cong [\K^{op}, \V](\phi, \A(S-,a))$,
\end{fact}
so that $\phi * S$ is equally the $\phi$-weighted 
limit of $S^{op}: \K^{op} \rightarrow \A^{op}$. Of course 
the limit $\{ \phi, T \}$ consists not just of the object
$\{ \phi, T \}$ but also of the representation \ref{fac2.1},
or equally of the corresponding {\em counit}
$\mu: \phi \rightarrow \A(\{\phi,T  \},T-)$; it is by 
{\em abus de langage}
that we usually mention only 
$\{\phi, T  \}$. When $\V = \Set$, we refind the classical 
(or ``conical'') limit of $T:\K^{op} \rightarrow \A$ and the classical
colimit of $S:\K \rightarrow \A$ as 
\begin{fact}\label{fac2.3}\spa
$lim\;T = \{ \Delta 1 , T \}$ $\;$ {\em and} $\;$ $colim\;S = \Delta 1 * S$
\end{fact} 
where $\Delta 1: \K^{op} \rightarrow \Set$ is the constant functor
at the one point set $1$.
Recall too that the weighted limits and colimits can be calculated 
using the classical ones when $\V = \Set$: for then the 
presheaf $\phi: \K^{op} \rightarrow \Set$ gives the discrete
op-fibration $d: el(\phi) \rightarrow \K^{op}$ where
$el(\phi)$ is the category of elements of $\phi$, and now
\begin{fact}\label{fac2.4}\spa
$\{ \phi, T \} = lim \{ \xymatrix{ el(\phi) \ar[r]^{d}&  \K^{op} \ar[r]^{T}
  & \A} \}$,
\end{fact}
\begin{fact}\label{fac2.5}\spa
$ \phi * S = colim \{ \xymatrix{ {el(\phi)}^{op} \ar[r]^{d^{op}}&  
\K \ar[r]^{S} & \A}    \}.$
\end{fact}
Recall finally that a functor $F: \A \rightarrow \B$ is said to {\em preserve
the limit} $\{ \phi, T \}$ as in \ref{fac2.1} when
$F(\{ \phi, T \})$ is the limit of $FT$ weighted by 
$\phi$, with counit 
$$\xymatrix{ \phi \ar[r]^(.3){\mu}
&  \A(\{ \phi, T \},T-) \ar[r]^(.45){F} 
& \B(F\{ \phi, T \},FT-) };$$ and $F$ is said to 
preserve the colimit
$\phi * S$ as in \ref{fac2.2} when $F^{op}$ preserves
$\{ \phi, S^{op} \}$.\\

We spoke above of a ``class $\Phi$ of colimits'' or a 
``class $\Psi$ of limits''; but this is loose and rather
dangerous language - the only thing that one can sensibly speak
of is {\em a class $\Phi$ of weights}. Then a category $\A$ {\em admits
$\Phi$-limits}, or is {\em $\Phi$-complete}, if $\A$ admits the limit
$\{\phi, T \}$ for each weight $\phi: \K^{op} \rightarrow \V$ in
$\Phi$ and each $T: \K^{op} \rightarrow \A$; while $\A$ {\em admits 
$\Phi$-colimits}, or is {\em $\Phi$-cocomplete}, when $\A$ admits 
the colimit $\phi * S$ for each $\phi: \K^{op} \rightarrow \V$
in $\Phi$ and each $S: \K \rightarrow \A$
(and thus when $\A^{op}$ is $\Phi$-complete). Moreover a 
functor $\A \rightarrow \B$ between $\Phi$-complete categories
is said to be {\em $\Phi$-continuous} when it preserves all 
$\Phi$-limits, and one defines {\em $\Phi$-cocontinuous} dually.
We write $\PhiCont$ for the 2-category of $\Phi$-complete
categories, $\Phi$-continuous functors, and all natural
transformations -- which is a (non full) sub-2-category
of $\VCAT$; and similarly $\PhiCo$ for the 2-category
of $\Phi$-cocomplete categories, $\Phi$-cocontinuous functors, 
and all natural transformations.\\
 
To give a class $\Phi$ of weights is to give, for each small
$\K$, those $\phi \in \Phi$ with domain $\K^{op}$; let us use
as in \cite{AK88} the notation
\begin{fact}\label{fac2.6}\spa
$\Phi[\K] = \{ \phi \in \Phi \mid dom (\phi) = \K^{op} \}$,
\end{fact}
so that
\begin{fact}\label{fac2.7}\spa
$\Phi = \Sigma_{\K\;small} \Phi[\K]$.
\end{fact}
In future, we look on $\Phi[\K]$ as a full subcategory
of the functor category $[\K^{op},\V]$ (which we may also 
call a {\em presheaf} category). The smallest class
of weights is the empty class $0$, and $\ZCont$ is just 
$\VCAT$. The largest class of weights consists of {\em all}
weights - that is, all presheaves with small domains - and 
we denote this class
by $\p$; the 2-category $\PCont$ is just the 2-category
$\Cont$ of {\em complete} categories and {\em continuous} functors,
and similarly $\PCo = \Co$.\\ 

There may well be different classes $\Phi$ and $\Psi$
for which the sub-2-categories $\PhiCont$ and $\PsiCont$
of $\VCAT$ coincide; which is equally to say that 
$\PhiCo$ and $\PsiCo$ coincide.
When $\V = \Set$, for instance, $\Cont = \PCont$ coincides
with $\PhiCont$ where $\Phi$ consists of the weights
for products and for equalizers. We define the 
{\em saturation} $\Phi^*$ of a class $\Phi$ of weights
as follows: the weight $\psi$ belongs to $\Phi^{*}$ when
every $\Phi$-complete category is also $\psi$-complete
and every $\Phi$-continuous functor is also 
$\psi$-continuous. 
Note that $\Phi^*$ was called in \cite{AK88} 
the {\em closure} of $\Phi$;
we now prefer the term ``saturation'', since ``closure''
already has so many meanings.
Clearly then, we have
\begin{fact}\label{fac2.8}\spa
$\PhiCont = \PsiCont\;\;\Leftrightarrow\;\;\PhiCo = \PsiCo\;\; 
\Leftrightarrow\;\;\Phi^* = \Psi^*.$
\end{fact} 
When $\V = \Set$, we can of course consider 
$\PhiCont$ where $\Phi$ consists of the $\Delta 1 : \K^{op} \rightarrow \Set$
for all $\K$ in some class $\D$ of small categories;
and we might write $\DCont$ for this 2-category $\PhiCont$
of {\em $\D$-complete} categories, {\em $\D$-continuous} functors, and 
all natural transformations. 
We underline the fact, however, that when $\V = \Set$,
NOT every $\PhiCont$ is of the form $\DCont$ for some $\D$ as 
above; a simple example of this situation is given in \cite{AK88}.\\

We spoke of $\VCAT$ as a 2-category, the category $\VCAT(\A,\B)$ 
having as its 
objects the $\V$-functors $T:\A \rightarrow \B$ and as its 
arrows the $\V$-natural transformations $\alpha: T \rightarrow S: \A
\rightarrow \B$. When $\A$ is small, however, we also have the 
$\V$-category $[\A,\B]$, whose underlying ordinary category $[\A,\B]_0$
is $\VCAT(\A,\B)$; an example is of course the presheaf 
$\V$-category $[\K^{op},\V]$ of \ref{fac2.1}.\\

When $\A$ is not small, $[\A,\B]$ may not exist as a $\V$-category, 
since the end $\int_a \B(Fa,Ga)$ giving the $\V$-valued hom 
$[\A,\B](F,G)$ may not exist in $\V$ for all $F,G: \A \rightarrow \B$. 
However it may exist for {\em some} pairs $F,G$, and then 
we can speak of $[\A,\B](F,G)$. This allows us the 
convenience of speaking of the 
limit $\{ \phi, T \}$ of \ref{fac2.1} or the colimit $\phi * S$ of 
\ref{fac2.2} even when $\K$ is not small (so that $\phi$ is 
no longer a weight, in the sense of this article) : for instance, 
we say that $\phi * S$ {\em exists} if the right side of \ref{fac2.2} 
{\em exists in $\V$ for each} $a$, {\em and is representable} as the left 
side of \ref{fac2.2}. In particular, we can speak, even when 
$\A$ is not small, of the 
possible existence of the left Kan extension $Lan_KT$ 
of some $T: \A \rightarrow \B$ along 
some $K: \A \rightarrow \C$, recalling from Chapter 4 
of \cite{Kel82} that it is given by
\begin{fact}\label{fac2.9}\spa
$Lan_KT(c)  \cong \C(K-,c)*T$,
\end{fact}
existing when the colimit on the right exists for each $c$.
\end{section}  

\begin{section}{Revision of the free $\Phi$-cocompletion of a category and
of saturated classes of weights}\label{section3} 
Another piece of background 
knowledge that we need to recall concerns  the ``left bi-adjoint'' to the 
forgetful 2-functor  $U_{\Phi}: \PhiCo \rightarrow \VCAT$.
(Note that it is convenient to deal with colimits
rather than limits.)\\

Recall from Section 4.8 of \cite{Kel82} that a presheaf 
$F:\A^{op} \rightarrow \V$, where $\A$ need not be small, 
is said to be {\em accessible} if it is the left Kan extension 
of some $\phi: \K^{op} \rightarrow \V$ with $\K$ small along 
some $H^{op}: \K^{op} \rightarrow \A^{op}$; which is to 
say, by \ref{fac2.9}, that $F$ has the form 
\begin{fact}\label{fac3.1}\spa
$Fa \cong \A(a,H-)* \phi $,
\end{fact}
which by (3.9) of \cite{Kel82} may equally be written as
\begin{fact}\label{fac3.2}\spa
$Fa \cong \phi* \A(a,H-) $.
\end{fact}
It is shown in Proposition 4.83 of \cite{Kel82} that, whenever 
$F$ is accessible, it is a left Kan extension as above for 
some $H$ that is {\em fully faithful}; in other words, that
$F$ is the left Kan extension of its restriction to some small 
full subcategory of $\A^{op}$.
Whenever $F$ is accessible, $[\A^{op},\V](F,G) $ exists for each $G$, an
easy calculation using the Yoneda isomorphism giving 
\begin{fact}\label{fac3.3}\spa
$\int_a [Fa,Ga] \cong \int_a [\phi*\A(a,H-), Ga] \cong [\K^{op},\V](\phi, GH^{op})$.
\end{fact}
Accordingly, for any $\A$ there is a $\V$-category 
$\p \A$ having as its objects the accessible
presheaves, and with its $\V$-valued hom given by the usual 
formula $\int_a [Fa,Ga]$; it was first introduced by Lindter 
\cite{Lin74}. Any presheaf $\A^{op} \rightarrow \V$
is accessible if $\A$ is small, 
being the left Kan extension of itself along the identity, so that 
$\p \A$ coincides with $[\A^{op},\V]$ for a small $\A$.\\

Every representable $\A(-,b)$ is accessible; we express it in the 
form \ref{fac3.2} by taking $\phi = I : \I^{op} \rightarrow \V$ and 
$H = b :\I \rightarrow \A$, where $\I$ is the unit $\V$-category and
$I$ is the unit for $\otimes$. Accordingly we have the fully-faithful
{\em Yoneda embedding}
$Y:\A \rightarrow \p \A$ sending $b$ to $\A(-,b)$, which
we sometimes loosely treat as an inclusion. Now 
calculating \ref{fac3.3} with $F = Yb$ gives at once the
{\em Yoneda isomorphism}
\begin{fact}\label{Yo} 
$\p \A(Yb, G) \cong Gb$.
\end{fact} 

By Proposition 5.34 of \cite{Kel82} the category $\p \A$ admits
all small colimits, these being formed pointwise from those in $\V$.
So the typical object $F$ of $\p \A$ as in \ref{fac3.2} can be written as
\begin{fact}\label{fac3.4}\spa
$F \cong \phi*YH$,
\end{fact}
this now being a colimit in $\p \A$. We can see \ref{fac3.4} as 
expressing the general accessible $F$ as a small colimit in $\p \A$
of representables.\\  

Recall from \cite{Kel82} p.154 that, given
a class $\Phi$ of weights and
a full subcategory $\A$ of a $\Phi$-cocomplete category 
$\B$, the {\em closure of $\A$ in $\B$
under $\Phi$-colimits} is the smallest full replete 
subcategory of $\B$ containing $\A$ and closed 
under the formation of $\Phi$-colimits in $\B$ -- namely the
intersection of all such.
For any class $\Phi$ of weights, and any category $\A$,
we write $\Phi(\A)$ for the closure of $\A$ in $\p A$
under $\Phi$-colimits,
with $Z:\A \rightarrow \Phi(\A)$ and 
$W: \Phi(\A) \rightarrow \p A$ for the full inclusions,
so that $Y:\A \rightarrow \p \A$ is the composite 
$WZ$; note that $W$ is $\Phi$-cocontinuous. We now 
reproduce (the main point of) 
\cite{Kel82} Theorem 5.35. The proof below is a 
little more direct than that given there, which referred
back to earlier propositions. The result itself must be older still, at
least for certain classes $\Phi$.

\begin{proposition}
\label{Pro3.5}
For any $\Phi$-cocomplete category $\B$,
composition with $Z$ gives an equivalence 
of categories
$$\PhiCo(\Phi(\A),\B) \rightarrow \VCat(\A,\B)$$
with an equivalence inverse given by the left Kan extension
along Z.
Thus $\Phi(-)$ provides a left bi-adjoint to the 
forgetful 2-functor $\PhiCo \rightarrow \VCAT$. 
\end{proposition}

\pf
By \ref{fac2.9}, the left Kan extension 
$Lan_ZG$ of $G:\A \rightarrow \B$ is given by 
$Lan_ZG(F) = \Phi(\A)(Z-,F)*G$, existing when this last 
colimit exists for each $F$ in $\Phi(\A)$. However
$\Phi(\A)(Z-,F) 
= \p \A(WZ-,WF) = \p \A(Y-,WF)$,
which by Yoneda is isomorphic to $WF$. 
Consider the full subcategory of $\Phi(\A)$
given by those $F$ for which $WF*G$ {\em does} exist; it 
contains the representables by Yoneda, and
it is closed under $\Phi$-colimits: for these exist in $\Phi(\A)$ 
and are preserved by $W$, while \cite{Kel82} (3.23) gives 
$(\phi*S)*G \cong \phi*(S-*G)$, either side existing if the 
other does; so the subcategory in question is all of $\Phi(\A)$.\\

What is more: $Lan_ZG = W-*G: \Phi(\A) \rightarrow \B$
preserves $\Phi$-colimits since $W$ does so and colimits
are cocontinuous in their weights (see \cite{Kel82} (3.23) again).
So one does indeed have a functor 
$Lan_Z: \VCAT(\A,\B) \rightarrow \PhiCo(\Phi(\A),\B)$, while one
has trivially the restriction functor 
$\PhiCo(\Phi(\A),\B) \rightarrow \VCAT(\A,\B)$ given 
by composition with $Z$.
By \cite{Kel82} (4.23), the canonical $G \rightarrow Lan_Z(G)Z$ 
is invertible for all $G$ since $Z$ is fully faithful. Thus 
it remains to consider the canonical
$\alpha: Lan_Z(SZ) \rightarrow S$ for a $\Phi$-cocontinuous 
$S: \Phi(\A) \rightarrow \B$. The $F$-component of
$\alpha$ for $F \in \Phi(\A)$ is the canonical 
$\alpha_F: WF * SZ \rightarrow SF$; and clearly
the collection of those $F$ for which $\alpha_F$ is 
invertible contains the representables and is closed 
under $\Phi$-colimits: therefore it is the totality
of $\Phi(\A)$.
\epf

\begin{remarks}
\label{Rem3.6}
We may express this by saying that 
$\Phi(\A)$ is the {\em free $\Phi$-cocomplete category on $\A$}.
As a particular case, $\p \A$ itself is the free 
cocomplete category on $\A$; in other words $\Phi(\A) = \p \A$ 
when $\Phi$ is the class of {\em all} weights - which is why
(identifying $\p(\A)$ with $\p \A$) we use $\p$ as the name 
for this class of all weights.
\end{remarks}

As shown in \cite{Kel82}, one can form $\Phi(\A)$ by transfinite induction. 
Define successively full replete subcategories $\A_{\alpha}$ of $\p \A$
as $\alpha$ runs through the ordinals: $\A_0$, which is equivalent 
to $\A$, consists of the representables, now in the sense of those 
presheaves isomorphic to some $\A(-,a)$; then $\A_{\alpha+1}$
consists of $\A_{\alpha}$ together with all $\Phi$-colimits
in $\p \A$ of diagrams in $\A_{\alpha}$; and for a limit ordinal
$\alpha$ we set $\A_{\alpha} = \bigcup_{\beta < \alpha} \A_{\beta}$.
This sequence stabilizes if, as we suppose, there exist
arbitrarily large inaccessible cardinals: for 
we have  $\Phi(\A) = \Phi_{\alpha}(\A)$
when $\alpha$ is the smallest regular cardinal
greater than $card(ob (\K))$ for all small $\K$ 
with $\Phi[\K]$ non-empty.
It follows that {\em $\Phi(\A)$ is a small category when $\A$  
and $\Phi$ are small}. 
In a number of important cases, one has $\Phi(\A) = \A_1$ in the 
notation above; it is so when $\Phi = \p$ since by \ref{fac3.4}
every accessible $F$ is a small colimit of representables,
and in the case $\V = \Set$ it is so by $\cite{Kel82}$
Theorem 5.37 when $\Phi$ consists of the weights for 
finite conical colimits.
However there is no special value in this condition,
which (as we shall see in \ref{3.13} below) always holds for a small $\A$ when 
the class $\Phi$ is {\em saturated}.\\ 

An explicit description of the saturation $\Phi^*$ of a class
$\Phi$ of weights was given by Albert and Kelly in \cite{AK88},
in the following terms:
\begin{proposition}
\label{Pro3.7} The weight $\psi: \K^{op} \rightarrow \V$ lies in the
saturation $\Phi^*$ of the class $\Phi$ if and only if the 
object $\psi$ of $\p \K = [\K^{op},\V]$ lies in the full 
subcategory $\Phi(\K)$ of $[\K^{op}, \V]$.
\end{proposition}
There is another useful way of putting this. When $\K$ is small,
both $\Phi[\K]$ and $\Phi(\K)$ make sense for any class 
$\Phi$; and in fact we have
\begin{fact}
\label{F3.8}\spa $\Phi[\K] \subset \Phi(\K),$
\end{fact}
since for $\phi: \K^{op} \rightarrow \V$ the Yoneda 
isomorphism 
\begin{fact}
\label{F3.9}\spa
$\phi \cong \phi * Y$
\end{fact} exhibits $\phi$ as 
an object of $\Phi(\K)$ when $\phi \in \Phi$.
We can write Proposition \ref{Pro3.7} as
\begin{fact} 
\label{F3.10}\spa
$\Phi^*[\K] = \Phi(\K)$,
\end{fact}
so that $\Phi$ is a saturated class precisely 
when 
\begin{fact} 
\label{F3.11}\spa
$\Phi[\K] = \Phi(\K)$
\end{fact}
for each small $\K$. In other words the class $\Phi$ is 
saturated precisely when, for each small $\K$, 
{\em the full subcategory
$\Phi[\K]$ of $[\K^{op},\V]$ contains the representables 
$\K(-,k)$ and is closed in $[\K^{op}, \V]$ under $\Phi$-colimits.}

\begin{example}\label{exsat}\end{example}
Consider the case when $\V$ is locally finitely presentable
as a closed category in the sense of \cite{Kel82-2},
and $\Phi$ is the class of finite weights as described
there; this includes the case where $\V= \Set$ and $\Phi$
is the set of weights for the classical finite colimits.
Then $\Phi^*[\K] = \Phi^*(\K)$ is the closure of $\K$
in $[\K^{op},\V]$ under finite colimits, which by \cite{Kel82-2}
Theorem 7.2 is the full subcategory of $[\K^{op},\V]$
given by the finitely presentable objects.\\ 

It follows of course from the definitions of $\Phi(\A)$
and of $\Phi^*$ that 
\begin{fact}\label{F3.12}\spa
$\Phi^*(\A) = \Phi(\A)$ 
\end{fact}
for any $\A$. We cannot write \ref{F3.11}
when $\K$ is replaced by a non-small $\A$, since then $\Phi[\A]$
has no meaning; but a partial replacement for it is provided
by the following, which was Proposition 7.4 in \cite{AK88}:
\begin{proposition}\label{3.13}
If the presheaf $F:\A^{op} \rightarrow \V$ lies in $\Phi(\A)$ for 
some saturated class $\Phi$, then $F$ is a $\Phi$-colimit 
in $\p A$ of 
representables; that is, $F \cong \phi * YH$ for some
$\phi: \K^{op} \rightarrow \V$ in $\Phi$ and some
$H: \K \rightarrow \A$. Since $\Phi$-colimits are 
formed in $\Phi(\A)$ as in $\p \A$, $F$ is equally
the colimit $\phi * ZH$ in $\Phi(\A)$.   
\end{proposition}
In other words an $F$ in $\Phi(\A)$ has the form 
\ref{fac3.4} with $\phi$ in $\Phi$. Equally,
this asserts that $F: \A^{op} \rightarrow \V$ is
the left Kan extension of $\phi: \K^{op} \rightarrow \V$
along $H^{op}: \K^{op} \rightarrow \A^{op}$. In fact,
we can take $H$ here to be fully faithful, as was shown
\cite{Kel82} Proposition 4.83 for the case $\Phi = \p$:
\begin{proposition}\label{F3.14}
For a saturated class $\Phi$, any $F$ in $\Phi(\A)$ is of the 
form $Lan_{H^{op}}\phi$ for some 
$\phi: \K^{op} \rightarrow \V$ in $\Phi$ and some 
fully faithful $H: \K \rightarrow \A$. 
\end{proposition}
\pf
We already have that $F \cong Lan_{T^{op}}\psi$ for
some $\psi: \LL^{op} \rightarrow \V$ in $\Phi$ and
some $T: \LL \rightarrow \A$. Let $T$ factorize as 
$T = HP$ where $H: \K \rightarrow \A$ is fully faithful
and $P: \LL \rightarrow \K$ is bijective on objects.
Then $\K$ is small since $\LL$ is small.
Now $F \cong Lan_{T^{op}}\psi \cong Lan_{H^{op}}\phi$,
where $\phi = Lan_{P^{op}}\psi$. However 
$\phi = Lan_{P^{op}}\psi \cong \psi * YP$, 
which, as a $\Phi$-colimit of representables,
lies in $\Phi(\K)$, and hence in $\Phi[\K]$.
\epf

It may be useful to understand extreme special cases
of one's notation.
First observe that the saturation $0^*$ of the empty class $0$ 
consists precisely of the representables -- that 
is, $0^*[\K] =0^*(\K)$ consists of the isomorphs
of the various  $\K(-,k): \K^{op} \rightarrow \V$. 
Another extreme case involves the empty $\V$-category $0$
with no objects. Of course $\p 0 = [0^{op},\V]$ is the 
terminal category $1$; its unique object is the unique 
functor $!:0^{op} \rightarrow \V$ and $1(!,!)$ is 
the terminal object $1$ of $\V$. (This differs in general
from the {\em unit $\V$-category} $\I$, with one object $*$
but with $\I(*,*)$$=$$I$.) 
So for any class $\Phi$, we have $\Phi[0] = 0$ if 
$!: 0^{op} \rightarrow \V$ is not in $\Phi$, and 
$\Phi[0] = [0^{op},\V] = 1$ otherwise. Now $\Phi(0)$ is the closure
of $0$ in $\p 0$ under $\Phi$-colimits, and any diagram
$T: \K \rightarrow 0$ has $\K=0$, so that $\Phi(0) = 0$
if $! \not \in \Phi$ and otherwise $\Phi(0)$ contains
$ ! * Y = !$, giving $\Phi(0) = 1$. So in fact 
$\Phi(0) = \Phi[0]$, being $0$ or $1$. Both are possible for a 
saturated $\Phi$; for $\p 0=1$, while the Albert-Kelly theorem
(Proposition \ref{Pro3.7} above) gives $0^*[0] = 0(0)= 0[0] = 0$.\\

Before ending this section, we recall a result characterizing 
$\Phi$-cocomplete categories, along with a short proof. This 
was Proposition 4.5 in \cite{AK88}.
\begin{proposition}
For any class $\Phi$ of weights, a category $\A$ admits 
$\Phi$-colimits if and only if the fully faithful
embedding $Z:\A \rightarrow \Phi(\A)$ admits a left adjoint;
that is, if and only if the full subcategory $\A$ given by
the representables is reflective in $\Phi(\A)$.
\end{proposition}
\pf
If $\A$ is reflective, it admits $\Phi$-colimits
because $\Phi(\A)$ does so. Suppose conversely that 
$\A$ admits $\Phi$-colimits, and write $\B$ for the full
subcategory of $\p \A$ given by those objects
admitting a reflection into $\A$; then $\B$ contains $\A$
and $\B$ is closed in $\p \A$ under $\Phi$-colimits
since $\A$ admits these; so that $\B$ contains $\Phi(\A)$,
as desired.
\epf
\end{section}

\begin{section}{Recognition theorems}\label{section4} 
We recall from Proposition 5.62 of \cite{Kel82}
a result characterizing categories of the form $\Phi(\A)$
-- or more precisely functors of the form 
$Z: \A \rightarrow \Phi(\A)$. At the same time, we give 
a direct proof; for the proof
in \cite{Kel82} refers back to earlier results in that book.

We begin with a piece of notation: for a category $\A$ and a 
class $\Phi$ of weights, we write $\A_{\Phi}$ for the full 
subcategory of $\A$ given by those $a \in \A$ for which the
representable $\A(a,-): \A \rightarrow \V$ preserves all
$\Phi$-colimits (That is, all $\Phi$-colimits that {\em exist in}
$\A$). There is no agreed name for $\A_{\Phi}$; the objects of
$\A_{\Phi}$ are usually called {\em finitely presentable}
when the $\Phi$-colimits are the classical filtered colimits;
while when $\Phi$ is the class $\p$ of all weights, the objects 
of $\A_\Phi$ were called {\em small projectives} in \cite{Kel82},
but have also been called {\em atoms} by some authors.
Let us use the name {\em $\Phi$-atoms} for the objects of $\A_{\Phi}$.
When $\A$ admits $\Phi$-colimits and hence $\Phi^*$-colimits,
it follows from the definition of $\A_{\Phi}$ that
\begin{fact}\label{fac4.1}\spa
$\A_{\Phi^*} = \A_{\Phi}.\label{F4.1}$ 
\end{fact}
For a functor $G: \A \rightarrow \B$, we get for each $b \in \B$
the presheaf $\B(G-,b): \A^{op} \rightarrow \V$. If this is accessible
for every $b$, we have a functor $\tilde{G}: \B \rightarrow \p \A$
in the notation of \cite{Kel82}. Note that $\tilde{G}G\cong Y:
\A^{op} \rightarrow \p \A$ when $G$ is fully faithful.

The following is the characterization result of Proposition
$5.62$ of \cite{Kel82} with a slightly expanded form of its 
statement.
\begin{proposition}\label{Pro4.1}
In order that $G:\A \rightarrow \B$ be equivalent to the free
$\Phi$-cocompletion $Z:\A \rightarrow \Phi(\A)$ of $\A$ for a
class $\Phi$ of weights, the following conditions are 
necessary and sufficient:\\
$(i)$ G is fully faithful (allowing us to treat $\A$
henceforth as a full subcategory of $\B$);\\ 
$(ii)$ $\B$ is $\Phi$-cocomplete;\\
$(iii)$ the closure of $\A$ in $\B$ under $\Phi$-colimits
is $\B$ itself;\\
$(iv)$ $\A$ is contained in the full subcategory $\B_{\Phi}$
of $\B$.\\
When these conditions hold, 
each functor $\B(G-,b): \A^{op} \rightarrow \V$ is accessible
and in fact lies in the full subcategory
$\Phi(\A)$ of $\p \A$.
Thus we have a functor $\tilde{G}: \B \rightarrow \p \A$
given by $\tilde{G}(b) = \B(G-,b)$, and this factorizes
as $WK$ where $W$ is, as before, the inclusion
from $\Phi(\A)$ to $\p \A$.
The functor $K: \B \rightarrow \Phi(\A)$ here is an equivalence, 
an equivalence inverse being given by 
$Lan_Z G : \Phi(\A) \rightarrow \B$, 
which by Proposition \ref{Pro3.5} is the unique 
$\Phi$-cocontinuous extension of $G$ to $\Phi(\A)$.
\end{proposition}
\pf The necessity of the first two conditions is clear. That
of the third results from the fact that the inclusion
$W: \Phi(\A) \rightarrow \p \A$ preserves $\Phi$-colimits 
by definition. 
For that of the fourth condition, 
the point is that $\Phi(\A)(Za,-) \cong \p \A(Ya,W-)$ preserves
$\Phi$-colimits: for $W$ does so, while
$\p \A (Ya,-): \p \A \rightarrow \V$ preserves
all small colimits, being isomorphic by Yoneda to the 
evaluation $E_a$.

We turn now to the proof of sufficiency.
First, to see that each $\B(G-,b)$ lies in the 
full replete subcategory $\Phi(\A)$ of $\p \A$, 
consider the full subcategory of $\B$ given by those 
$b$ for which this is so; this contains $\A$ since 
$\B(G-,Ga) \cong Ya$ because $G$ is fully faithful, 
and it is closed in $\B$ under $\Phi$-colimits 
by $(iv)$, since $\Phi(\A)$ is closed 
under these in $\p \A$; so it is all of $\B$. 

Thus we have indeed a functor $K: \B \rightarrow \Phi(\A)$,
sending $b$ to $\B(G-,b)$. We next show that $K$
or equivalently $\tilde{G} = WK: \B \rightarrow \p \A$ 
is fully faithful.
Consider the full subcategory of $\B$ given by those $b$ 
for which the map 
$\tilde{G}_{b,c}: \B(b,c) \rightarrow \p \A(\tilde{G}(b),\tilde{G}(c))$ 
is invertible for all $c$. 
We observe that it contains $\A$ since $G$ is fully faithful,
and that it is closed under $\Phi$-colimits since $\A \subset \B_{\Phi}$.
Thus it is all of $\B$.

It remains to show that $K$ and 
$S = Lan_{Z}G: \Phi(\A) \rightarrow \B$ are equivalence-inverses.
Recall that $\tilde{G}G \cong Y$ since $G$ 
is fully faithful.
Also recall from Proposition \ref{Pro3.5} that $S$ is the
essentially unique $\Phi$-cocontinuous functor with $SZ \cong G$.
So $WKSZ \cong \tilde{G}G \cong Y \cong WZ$, giving $KSZ \cong Z$ 
since $W$ is fully faithful, and then giving $KS \cong 1$ by Proposition 
\ref{Pro3.5} since 
$KS$ and $1$ are $\Phi$-cocontinuous, $K$ being so because 
$\A \subset \B_{\Phi}$.
Finally $KS \cong 1$ gives $KSK \cong K$; whence $SK \cong 1$
since $K$ (as we saw) is fully faithful.
\epf

Proposition \ref{Pro4.1} is of particular interest
in the case of a small $\A$. We may cast the result for a small 
$\A$ in the form: 
\begin{proposition}\label{Pro4.2}
For a class $\Phi$ of weights, the following properties
of a category $\B$ are equivalent:\\
$(i)$ For some small $\K$, there is an equivalence 
$\B \simeq \Phi(\K)$;\\
$(ii)$ $\B$ is $\Phi$-cocomplete and has a small full subcategory
$\A \subset \B_{\Phi}$ such that every object of $\B$ is a 
$\Phi^*$-colimit of a diagram in $\A$;\\
$(iii)$ $\B$ is $\Phi$-cocomplete and has a small full
sub-category $\A \subset \B_{\Phi}$ such that the closure
of $\A$ in $\B$ under $\Phi$-colimits is $\B$ itself.\\
Under the hypothesis $(iii)$ - and so a fortiori under $(ii)$ -
if $G: \A \rightarrow \B$ denotes the inclusion, the functor
$\tilde{G}: \B \rightarrow \p \A$ is fully faithful, with 
$\Phi(\A)$ for its replete image.
\end{proposition}

\begin{remarks}\label{Rem4.3}\end{remarks}
(a) When $\Phi$ is the class $\p$ of all weights, we get a
characterization here of the functor category $\p \K = [\K^{op},\V]$
for a small $\K$; note that it differs from the characterization given
in \cite{Kel82} Theorem 5.26, which replaces the condition that 
$\B$ be the colimit closure of $\A$ by the condition that $\A$ be 
strongly generating in $\B$; but these conditions are very similar
in strength by \cite{Kel82} Proposition 3.40.\\

\noindent(b) Theorem 5.3 of \cite{BQR98} is the special case where $\V$ is 
locally finitely presentable as a closed category in the sense
of \cite{Kel82-2} and $\Phi$ is the saturated class of $\alpha$-flat 
presheaves. (See Section \ref{section5} below.)\\

Let us mention the following consequence of Proposition \ref{Pro4.2}.
\begin{proposition}\label{Pro3.15}
For a small $\K$ and a saturated class $\Phi$,
let $\A$ be a full reflective subcategory 
of $\Phi(\K)$ that is closed in $\Phi(\K)$ under $\Phi$-colimits.
Then $\A$ is equivalent to $\Phi(\LL)$ for a small $\LL$.
\end{proposition}
\pf
Write $J: \A \rightarrow \Phi(\K)$ for the inclusion,
with $R: \Phi(\K) \rightarrow \A$ for its left adjoint, 
and regard $Z: \K \rightarrow \Phi(\K)$ as an 
{\em inclusion} of the representables in $\Phi(\K)$.
The objects $RZk$ of $\A$ with $k \in \K$ constitute 
a full subcategory $\LL$ of $\A$. By hypothesis, $\A$ 
admits $\Phi$-colimits and $J$ preserves these.
The subcategory $\LL$ lies in $\A_{\Phi}$, 
because $\A(RZk,-) \cong \Phi(\K)(Zk,J-)$ preserves 
$\Phi$-colimits since both $J$ and 
$\Phi(\K)(Zk,-)$ (being the evaluation at $k$) do
so. Finally every object $a$ of $\A$ is a $\Phi$-colimit
of a diagram taking its values in $\LL$; for $Ja \in \Phi(\K)$
is a $\Phi$-colimit $Ja * Z$, and $R$ preserves this colimit,
so that $a \cong RJa \cong Ja* RZ$, where the diagram 
$RZ: \K \rightarrow \A$ takes its values in $\LL$.
\epf
\end{section}

\begin{section}{Limits and colimits commuting in $\V$}\label{section5}
The new observations to which we now turn begin with the general
study of the commutativity in $\V$ of limits and colimits.
For a pair of weights $\psi: \K^{op} \rightarrow \V$ and 
$\phi: \LL^{op} \rightarrow \V$, to say that
\begin{fact}\label{fac5.1}\spa 
$\phi * - : [\LL,\V] \rightarrow \V$ preserves $\psi$-limits
\end{fact}
is equally to say that 
\begin{fact}\label{fac5.2}\spa 
$\{ \psi, - \} : [\K^{op},\V] \rightarrow \V$ preserves $\phi$-colimits,
\end{fact}
because each in fact asserts the invertibility, for every functor 
$S: \K^{op} \otimes \LL \rightarrow \V$, of the canonical 
comparison morphism
\begin{fact}\label{fac5.3}\spa
$\phi? * \{ \psi-, S(-,?) \} \rightarrow \{ \psi-, \phi? * S(-,?) \}.$
\end{fact}

When these statements are true for every such $S$, we say
that {\em $\phi$-colimits commute with $\psi$-limits in $\V$}.
For classes $\Phi$ and $\Psi$ of weights, if \ref{fac5.1} (or 
equivalently \ref{fac5.2}) holds for all $\phi \in \Phi$ and all 
$\psi \in \Psi$, we say that 
{\em $\Phi$-colimits commute with $\Psi$-limits in $\V$}.
For any class $\Psi$ of weights we may 
consider the class $\Psi^+$ of all weights $\phi$ for which 
$\phi$-colimits commute with $\Psi$-limits in $\V$; and for any class
$\Phi$ of weights we may consider the class $\Phi^-$ of all weights
$\psi$ for which $\Phi$-colimits commute with $\psi$-limits
in $\V$. We have here of course a Galois connection,
with $\Phi \subset \Psi^+$ if and only if 
$\Psi \subset \Phi^-$.
Note that $[\LL, \V]$ and $\V$ in \ref{fac5.1} admit all 
(small) limits; so that by the definition above 
of the saturation $\Psi^*$ of a class $\Psi$ of weights,
if $\phi * -$ preserves all $\Psi$-limits, it also 
preserves $\Psi^{*}$-limits. From this and a dual 
argument, one concludes that:
\begin{proposition}\label{pro5.1}For any classes
$\Phi$ and $\Psi$ of weights, the classes $\Phi^-$ and 
$\Psi^+$ are saturated; so that 
$\Psi^{+\;*} = \Psi^+$ and $\Phi^{-\;*} = \Phi^-$.
Moreover $\Psi^+ = \Psi^{*\;+}$ and
$\Phi^- = \Phi^{*\;-}$.
\end{proposition} 

When $\Psi$ consists of the weights for finite limits
(in the usual sense for ordinary categories, or in the sense
of \cite{Kel82-2} when $\V$ is locally finitely presentable 
as a closed category), it has been customary to call the 
elements of $\Psi^+$ the {\em flat} weights, as they
are those $\phi$ having $\phi * -: [\LL,\V] \rightarrow \V$
left exact. (Note the corresponding use of ``$\alpha$-flat''
in Definition 4.1 of \cite{BQR98}.) Accordingly for a general 
$\Psi$ we call the elements of $\Psi^+$ the {\em $\Psi$-flat} 
weights.\\   

Recall that the limit functor $\{\psi, - \}$ of 
\ref{fac5.2} is just the representable functor
$[\K^{op},\V](\psi,-):[\K^{op},\V] \rightarrow \V$.
Accordingly $\psi:\K^{op} \rightarrow \V$ lies
in ${\Phi}^-$ for a given class $\Phi$ if and 
only if it lies in the subcategory $[\K^{op},\V]_{\Phi}$:
\begin{fact}\label{fac5.6}\spa
$\Phi^-[\K] = \Phi^-(\K) = [\K^{op},\V]_{\Phi}$. 
\end{fact}
In other words, the elements $\psi$ of $\Phi^-[\K]$ are
the {\em $\Phi$-atoms} of $[\K^{op},\V]$; we also call
them the {\em $\Phi$-atomic} weights. When $\Phi$ is the class
$\p$ of all weights, the elements of $\p^-$ are also
called the {\em small projective} weights.\\ 

Part of the saturatedness of ${\Phi}^-$ -- namely the 
closedness of ${\Phi}^-[\K]$ in $[\K^{op}, \V]$ under
$\Phi^-$-colimits -- is the special case for
$\A = [\K^{op},\V]$ of the following more general result:
\begin{proposition}\label{pro5.2} For any class $\Phi$ of weights 
and any category $\A$, the full subcategory $\A_{\Phi}$
of $\A$ is closed in $\A$ under any $\Phi^-$-colimits 
that exist in $\A$.
\end{proposition}
\pf Let the colimit $\psi * S$ exist, where 
$\psi:\K^{op} \rightarrow \V$ lies in $\Phi^-$ and
$S:\K \rightarrow \A$ takes its values in $\A_{\Phi}$.
Then by definition 
$$\A(\psi * S,a) \cong [\K^{op}, \V](\psi, \A(S-,a)).$$
Since each $\A(Sk,-)$ preserves $\Phi$-colimits,
and since $[\K^{op},\V](\psi,-)$ preserves $\Phi$-colimits
by \ref{fac5.6}, it follows that $\A(\psi*S,-)$ preserves
$\Phi$-colimits: that is to say $\psi*S \in \A_{\Phi}$. 
\epf

\begin{example}\label{exsat1}\end{example} 
When $\V = \Set$, let $\Psi$ be the class 
of weights for (classical conical) finite limits: that is,
the set of all $\Delta 1: \K^{op} \rightarrow \Set$ with 
$\K$ finite. Then $\Psi^+$ consists of those $\phi: \LL^{op} \rightarrow
\Set$ with $\phi * -: [\LL,\Set] \rightarrow \Set$ left exact;
that is, the {\em flat} presheaves 
$\phi: \LL^{op} \rightarrow \Set$. As is well known, these 
are those presheaves $\phi$ for which ${(el(\phi))}^{op}$
is filtered. 
Since $\phi*S$ for $S: \LL \rightarrow \A$ is given
as in \ref{fac2.5} by
$colim \{ 
\xymatrix{ {el(\phi)}^{op} \ar[r]^{ d^{op} } & \LL \ar[r]^S & \A } 
\}$, a functor $[\K^{op},\Set] \rightarrow \Set$
is $\Psi^+$-cocontinuous if and only if it preserves 
filtered colimits; that is, if and only if it is 
finitary. By \ref{fac5.6}, therefore, $\Psi^{+-}$ consists 
of those $\psi: \K^{op} \rightarrow \V$ for which 
$[\K^{op},\Set](\psi,-)$ preserved filtered colimits;
that is, those $\psi$ that are finitely presentable 
in $[\K^{op},\Set]$. It follows from \ref{exsat}
that $\Psi^{+-}$ coincides in this case with $\Psi^*$.

\begin{example}\label{exsat2}\end{example}
With $\V = \Set$ again, let $\Psi$ consist of the single object 
$0^{op} \rightarrow \Set$, 
where $0$ is the empty category: so a $\Psi$-limit is a terminal 
object. Now $\phi: \LL^{op} \rightarrow \Set$ lies in 
$\Psi^{+}$ whenever $\phi*-: [\LL, \Set] \rightarrow \Set$ preserves 
the terminal object; which is to say that $\phi * \Delta 1 \cong 1$, 
or equally that $colim(\phi) \cong 1$, or again that 
$el(\phi)$ is connected. So the presheaf 
$\psi: \K^{op} \rightarrow \Set$ lies in 
$\Psi^{+-}$ just when $[\K^{op}, \Set](\psi,-)$ preserves 
connected (conical) colimits. This time $\Psi^{+-}$ strictly 
includes $\Psi^*$. 
For $\Psi^*(\K)$, being the closure of the representables in 
$[\K^{op}, \Set]$ 
under $\Psi$-colimits, consists of the representables together with the 
initial object $\Delta 0: \K^{op} \rightarrow \Set$. 
When $\K$ has one object, being given 
by the monoid $\{1,e\}$ with $e^2 = e$, the subcategory 
$\Q(\K)$ of $[\K^{op},\Set]$ 
given by the Cauchy completion of $\K$ has, by Section 5.8 of 
\cite{Kel82}, two objects, the representable object $*$ and the 
equalizer $E$ of the two maps $1,e : * \rightarrow *$, which splits 
the idempotent $e$; and $E$ is not $\Delta 0$ since there is an arrow 
from $*$ to $E$ because $ee = e$. Now $\Q = \p^-$ by Section 6 below, 
and $\p^- \subset \Psi^{+-}$ because $\Psi^+ \subset \p$. 
So in this case, there are objects of $\Psi^{+-}(\K)$ 
which are not contained in $\Psi^*(\K)$, and 
$\Psi^*$ is properly contained in $\Psi^{+-}$.

\begin{remark} \end{remark}
When $\V = \Set$, it is well known (see for example Theorem 5.38
of \cite{Kel82}) that the flat weights $\K^{op} \rightarrow \V$
are precisely the filtered conical colimits of representables,
and hence constitute the closure of $\K$ in $[\K^{op},\V]$
under filtered conical colimits. This is false for a general
$\V$ that is locally finitely presentable as a closed category;
if \cite{BQR98} seems to suggest otherwise, it is only because
those authors {\em define} ``filtered colimit'' to mean
``colimit weighted by a flat weight''.

\end{section}

\begin{section}{The class $\Q$ of small projective weights}\label{section6}
This section is devoted to the study of the saturated class
$\Q = \p^-$ of {\em small projective} weights. 
So for a small $\K$, \ref{fac5.6} gives 
\begin{fact}\spa
$\Q[\K] = \Q(\K) = {[\K^{op}, \V]}_{\p}$, 
\end{fact}
consisting of those $\phi: \K^{op} \rightarrow \V$ for 
which $\{ \phi, - \} = [\K^{op},\V](\phi,-): [\K^{op},\V] \rightarrow \V$ 
preserves all small colimits. We shall establish the following 
alternative characterizations of $\Q$.
First from Proposition \ref{fladj2} below:
\begin{fact}\label{facadj}\spa
$\phi: \K^{op} \rightarrow \V$ lies in $\Q$ if and only 
the corresponding module $\xymatrix{ \I \ar[r]|{\circ} & \K}$
is a left adjoint.
\end{fact}
Proposition \ref{yac} below gives:
\begin{fact}\spa 
$\Q$ is the class $\p^+$ of $\p$-flat weights.
\end{fact}
Finally, as Street showed in \cite{Str83},
\begin{fact}\spa 
$\Q$ is the class of {\em absolute} weights.
\end{fact}
Moreover there is an adjunction 
$L \dashv R: [\K,\V]^{op} \rightarrow [\K^{op},\V]$, 
due in the case
$\V = \Set$ to Isbell, which restricts to an equivalence 
${( \Q(\K^{op}))}^{op} \simeq \Q(\K)$
between the full subcategories of small projectives in 
$[\K,\V]$ and in $[\K^{op},\V]$. 
In terms of modules, this equivalence sends a right adjoint 
module $\xymatrix{\K \ar[r]|{\circ} & \I}$ to its left adjoint 
$\xymatrix{\I \ar[r]|{\circ} & \K}$.\\ 

Recall that by a module $\xymatrix{\A \ar[r]|{\circ} &  \B}$
is meant a functor $\B^{op} \otimes \A \rightarrow \V$ with $\A$ and
$\B$ small,
and that modules with their usual composition and 2-cells
form a bicategory $\VMod$. Recall further that
each functor $T: \A \rightarrow \B$ gives 
rise to modules $T_*: \xymatrix{\A \ar[r]|{\circ} & \B}$
and $T^*: \xymatrix{\B \ar[r]|{\circ} & \A}$, where
\begin{fact}\label{M1.8}\spa
$T_*(b,a) = \B(b,Ta)$$\;\;\;$and$\;\;\;$$T^*(a,b) = \B(Ta,b),$
\end{fact}
and that $T_*$ is left adjoint to $T^*$ in $\VMod$.
Recall finally that the bicategory $\VMod$ is {\em closed},  
admitting all {\em right liftings} and all {\em right extensions} 
as follows:  given modules
$f: \xymatrix{\A \ar[r]|{\circ} & \B}$,
$g: \xymatrix{\C \ar[r]|{\circ} & \A}$
and $h: \xymatrix{\C \ar[r]|{\circ} & \B}$,
we have the right lifting 
$\Cl f,h \Cr: \xymatrix{\C \ar[r]|{\circ} & \A}$
of $h$ through $f$
and the right extension
$\Sl g,h \Sr: \xymatrix{\A \ar[r]|{\circ} & \B}$ 
of $h$ along $g$,
given by:
\begin{fact}\label{M1.2}\spa
$\Cl f,h \Cr(a,c) = \int_b [f(b,a),h(b,c)]$
\end{fact}
and 
\begin{fact}\label{M1.3}\spa
$\Sl g,h \Sr (b,a) = \int_c [g(a,c),h(b,c)]$,
\end{fact}
satisfying the universal properties
\begin{fact}\label{M1.4}\spa
$\VMod(\A,\B)(f,\Sl g,h \Sr) \cong
\VMod(\C,\B)(fg,h) \cong 
\VMod(\C,\A)(g,\Cl f,h \Cr).$
\end{fact}
The second isomorphism corresponds by Yoneda to a morphism 
$\epsilon: f \Cl f,h \Cr   \rightarrow h: \C \rightarrow \B$
which is said, in the language of \cite{StWa78},
to {\em  exhibit} $\Cl f,h \Cr$ 
{\em as the right lifting of $h$ through $f$}. 
Such a lifting  $\Cl f,h \Cr$ is {\em respected} by a 
$k: \xymatrix{ \D \ar[r]|{\circ} & \C }$
when the 2-cell $\epsilon k$ exhibits $\Cl f, h \Cr k$ as the 
right lifting $\Cl f, hk  \Cr$ of $h  k$ through $f$, and 
the lifting  $\Cl f,h \Cr$ is {\em absolute}
when it is respected by every such arrow $k$.\\

As in any closed bicategory, we have the following 
characterization of left adjoints:
\begin{proposition}\label{sladj1} In $\VMod$, the following 
statements are equivalent:\\
$(i)$ $f: \xymatrix{\A \ar[r]|{\circ} & \B}$ has a right adjoint;\\
$(ii)$ for all $h: \xymatrix{\C \ar[r]|{\circ} & \B}$, 
the right lifting $\Cl f,h \Cr$ of $h$ through $f$ is absolute;\\
$(iii)$ the right lifting $\Cl f,1 \Cr$ of 
$1: \xymatrix{\B \ar[r]|{\circ} & \B}$ through $f$ is 
respected by $f$.\\
When these are satisfied, the right adjoint $f^*$ of $f$ is
the right lifting $\Cl f,1 \Cr$ of $1$ through $f$;
moreover the right lifting $\Cl f,h \Cr$ of $(ii)$ above
is given by $f^* h$.
\end{proposition}
There exists of course a dual characterization of right adjoints, 
in terms of right extensions. Thus letting  
$h$ in $\ref{M1.4}$ be $1_{\B}: \xymatrix{ \B \ar[r]|{\circ} & \B}$
yields an adjunction
\begin{fact}\label{genadj}\spa
$\Cl  -,1 \Cr \dashv \Sl -, 1 \Sr : \VMod(\B,\A)^{op} \rightarrow 
{\VMod(\A,\B)},$ 
\end{fact}
which restricts to an equivalence between the right adjoints
$\xymatrix{ \B \ar[r]|{\circ} & \A}$ and the left adjoints 
$\xymatrix{ \A \ar[r]|{\circ} & \B}$, for 
$\Cl -,1 \Cr$ sends a left adjoint to its right adjoint,
while $\Sl -,1 \Sr$ sends a right adjoint 
to its left adjoint.\\

We now translate Proposition \ref{sladj1} into the language  
of functors.
Consider again morphisms 
$f: \xymatrix{ \A \ar[r]|{\circ} &  \B }$,
$h: \xymatrix{ \C \ar[r]|{\circ} & \B }$ and
$k: \xymatrix{ \D \ar[r]|{\circ} & \C }$ in $\VMod$.
These correspond respectively to functors
$F: \A \rightarrow [\B^{op},\V]$, $H: \C \rightarrow [\B^{op},\V]$,
and $K: \D \rightarrow [\C^{op},\V]$; let us also
write $H': \B^{op} \rightarrow [\C,\V]$ for the 
other functor corresponding to $h$.
One checks straightforwardly that
\begin{fact}\label{rlpre}\spa 
$k$ respects the right lifting $\Cl f,h \Cr$ of $h$ through $f$
\end{fact} 
is equivalent to
\begin{fact}\label{colprelim}\spa
for all $a$ in $\A$ and all $d$ in $\D$,
$Kd * - : [\C,\V] \rightarrow \V$
preserves the limit $\{ Fa , H' \}$,
\end{fact}
which, by the equivalence of \ref{fac5.1} and \ref{fac5.2},
is further equivalent to
\begin{fact}\label{limprecol}\spa 
for all $a$ in $\A$ and all $d$ in $\D$,
the colimit $Kd * H$ is preserved
by $\{ Fa, - \}: [{\B}^{op},\V] \rightarrow \V$.
\end{fact} 
Our particular interest is in the case $\A = \I$
of the above: to give a module 
$f: \xymatrix{ \I \ar[r]|{\circ} & \B}$ 
is to give a presheaf $\phi: \B^{op} \rightarrow \V$,
and we write $\overline{\phi}$ for $f$.
Equally to give a module $g: \xymatrix{ \B \ar[r]|{\circ} & \I}$
is to give a presheaf $\psi: \B \rightarrow \V$, and
we write $\underline{\psi}$ for $g$.
Now \ref{sladj1} gives the following proposition,
in which the assertion $(ii)$ is the direct translation  
of the fact that for any module 
$h: \xymatrix{\C \ar[r]|{\circ} & \B}$ 
the right lifting $\Cl f,h \Cr$
is respected by any module 
$\xymatrix{ \I \ar[r]|{\circ} & \C}$.
\begin{proposition}\label{fladj2}
Given a weight $\phi: {\B}^{op} \rightarrow \V$, the following
conditions are equivalent:\\
$(i)$ $\overline{\phi}$ has a right adjoint $\underline{\psi}$;\\
$(ii)$ the representable functor 
$\{\phi ,-\} = [{\B}^{op}, \V](\phi,-): [{\B}^{op},\V] \rightarrow \V$ is 
cocontinuous, that is, $\phi$ is a small projective;\\
$(iii)$ the representable
functor $\{\phi,-\}: [{\B}^{op}, \V] \rightarrow \V$
preserves the colimit $\phi * Y$ of $Y:\B \rightarrow [{\B}^{op},\V]$
weighted by $\phi:{\B}^{op} \rightarrow \V$.\\
When these are satisfied, a right adjoint $\underline{\psi}$
of $\overline{\phi}$ is given by taking  
\end{proposition}
\begin{fact}\label{B2.17bis}
$\;\;$ $\psi = [\B^{op},\V](\phi,Y-).$
\end{fact}
Dually, $\underline{\psi}$ has a left adjoint
if and only $\psi$ is a small projective in $[\B,\V]$,
and then a left adjoint $\overline{\phi}$ of $\underline{\psi}$
is given by 
\begin{fact}\label{B2.18bis}
$\;\;$$\phi = [\B,\V](\psi,Y'-)$,
\end{fact}
where $Y'$ is the Yoneda embedding $\B^{op} \rightarrow [\B,\V]$.\\

Recall that every functor $G:\B \rightarrow \C$ where
$\B$ is small and $\C$ is cocomplete has the essentially unique
cocontinuous extension $Lan_Y G = -*G:[\B^{op},\V] \rightarrow \C$ 
along the Yoneda embedding $Y: \B \rightarrow [\B^{op},\V]$,
and $-*G$ has in fact the right adjoint 
$\tilde{G}: \C \rightarrow [\B^{op},\V]$ given by
$\tilde{G}(c) = \C(G-,c)$.
Moreover $\tilde{G}G$ is isomorphic to $Y$ when 
$G$ is fully faithful.
Applying this when $G$ is the Yoneda embedding 
${Y'}^{op}: \B \rightarrow [\B,\V]^{op}$, we get a 
commutative diagram
$$\xymatrix{[\B,\V]^{op} \ar@<1ex>[rr]^{R} 
& 
& [\B^{op},\V] \ar@<1ex>[ll]^{L}\\
\\
& \B \ar[luu]^{{Y'}^{op}} \ar[ruu]_Y & 
}$$
with $L$ left adjoint to $R$;
an easy calculation gives 
\begin{fact}\label{colim}\spa
$L(\phi) = [\B^{op},\V](\phi,Y-)$ {\em and} 
$R(\psi) = [\B,\V](\psi,Y'-)$.
\end{fact}
This adjunction, which we shall call the {\em Isbell
adjunction}, is in fact the case $\A = \I$
of the adjunction \ref{genadj}. Moreover 
when $\phi \in [\B^{op},\V]$ is a small 
projective, it follows from \ref{B2.17bis} 
that $L(\phi) = \psi$  
where the module $\underline{\psi}$ is the right adjoint
of the module $\overline{\phi}$; so that in fact $\psi$
too is a small projective.
Dually, when $\psi \in [\B,\V]$ is a small projective, 
it follows from \ref{B2.18bis} that 
$R(\psi) = \phi$ where $\overline{\phi}$ is the left adjoint
of $\underline{\psi}$; with $\phi$ too a small projective. 
In other words the adjunction $L \dashv R$
restricts to an equivalence at the level of small projectives,
which we may write as 
\begin{fact}\label{protopro}
$ {( \Q(\B^{op}) )}^{op} \simeq {\Q( \B ) }.$
\end{fact}
If $\phi \in [\B^{op},\V]$ and $\psi \in [\B,\V]$
are small projectives which correspond in this
equivalence, the functor $[\B,\V](\psi,-): [\B,\V] \rightarrow \V$,
being cocontinuous, has the form $-*\theta$ where $\theta$ 
is its composite with $Y': \B^{op} \rightarrow [\B,\V]$.
However this composite is $\phi$ by $6.16$,
and we can write $- * \phi$ as $\phi * -$; so we have 
\begin{fact}\label{flem}\spa
$[\B,\V](\psi, -) \cong \phi * - : [\B,\V] \rightarrow \V$.
\end{fact}
In terms of modules, this is just the observation 
that a right extension along a right adjoint 
is given by composition with its left adjoint --
since for a $\phi$ and a $\psi$ as above we have 
an adjunction 
$\overline{\phi} \dashv \underline{\psi}: 
\xymatrix{\B \ar[r]|{\circ} & \I}$.
This leads to another characterization of the small 
projectives:
\begin{proposition}\label{yac}
For a weight $\phi: \B^{op} \rightarrow \V$, 
the following conditions are equivalent:\\
$(i)$ $\phi$ is a small projective;\\
$(ii)$ $\phi * - : [\B,\V] \rightarrow \V$ is representable;\\
$(iii)$ $\phi* - : [\B,\V] \rightarrow \V$ is continuous;\\
$(iv)$ $\phi * - : [\B,\V] \rightarrow \V$  
preserves the limit $\{ \phi , Y' \}$ of 
$Y': \B^{op} \rightarrow [\B,\V]$ weighted by $\phi$. 
\end{proposition}
\pf
$(i) \Rightarrow (ii)$ by \ref{flem}.
It is trivial that $(ii) \Rightarrow (iii) \Rightarrow (iv)$.
By the equivalence of \ref{fac5.1} and \ref{fac5.2}, 
$(iv)$ is equivalent to the preservation by $\{\phi, -\}$ of the
colimit $\phi * Y$; and this is equivalent to $(i)$ by 
Proposition \ref{fladj2}.
\epf

\begin{remark}\end{remark}\noindent
The assertion $(iii)$ of the proposition above
may be expressed by saying that $\Q$ is the class $\p^{+}$ 
of $\p$-flat weights.\\

There is a further characterization of the weights in $\Q$, 
due to Street. A weight $\phi: \B^{op}\rightarrow \V$ 
is said to be {\em absolute} if each limit $\{ \phi, T \}$, 
where $T: \B^{op} \rightarrow \C$ say, is preserved by {\em every} 
functor $P: \C \rightarrow \D$; or equally if each colimit
$\psi * S$, where $S: \B \rightarrow \C$, is preserved 
by every functor $P: \C \rightarrow \D$. Street showed the following,
in a context wider than ours, in \cite{Str83};
we give a proof (in our context) for completeness:
\begin{theorem}\label{Str}
A weight $\phi: \B^{op} \rightarrow \V$
is absolute precisely when it is a small projective 
in $[\B^{op},\V]$.
\end{theorem}
\pf
One direction is clear:
to say that $\phi: \B^{op} \rightarrow \V$ is a small
projective is, by the equivalence of \ref{fac5.1} and \ref{fac5.2},  
to say that, for each 
$T: \B^{op} \rightarrow [\A,\V]$ with $\A$ small
and for each weight $\psi: \A^{op} \rightarrow \V$, the limit
$\{ \phi, T \}$ is preserved by the functor 
$\psi * - : [\A,\V] \rightarrow \V$; so that each absolute
$\phi: \B^{op} \rightarrow \V$ is certainly a small projective
in $[\B^{op},\V]$.\\

As a preliminary to the proof of the converse,
recall that the defining property
of the colimit $\phi * S$ for $S: \B \rightarrow \C$
is an isomorphism
\begin{fact}\label{fac5.1sup}
$\;\;$$\C(\phi * S,c ) \cong [\B^{op},\V](\phi, \C(S-,c)).$
\end{fact}
However $\C(Sb,c) = S^*(b,c)$; and then if 
$\phi: {\B}^{op} \rightarrow \V$ corresponds to the module
$\overline{\phi}: \xymatrix{ \I \ar[r]|{\circ} &  \B}$, 
the right side of \ref{fac5.1sup} is 
$\Cl \overline{\phi}, S^* \Cr(*,c)$, where $*$ denotes the 
unique object of $\I$. Finally the object $\phi * S$ of $\C$ 
corresponds to a functor $\phi * S: \I \rightarrow \C$ and 
hence to a module 
${(\phi * S)}^*: \xymatrix{\C \ar[r]|{\circ} &  \I}$
with $(\phi * S)^*( * , c) = \C( \phi * S, c)$; so that the defining 
equation \ref{fac5.1sup} of $\phi*S$ may be written as
\begin{fact}\label{fac5.2sup}
$\;\;$${(\phi * S)}^* \cong \Cl \overline{\phi}, S^*  \Cr$
\end{fact}
which is just to say that the lifting 
of $S^*$ through $\overline{\phi}$ is given
by ${(\phi * S)}^*$. 

To ask $P: \C \rightarrow \D$ to preserve the  
colimit $\phi * S$ is to ask the invertibility of the canonical 
comparison $ \phi * (PS) \rightarrow P(\phi * S)$  or equally 
of the canonical comparison 
${(\phi * S)}^* P^* \rightarrow {(\phi * PS)}^*$. By
\ref{fac5.2sup} this may be written in the form
\begin{fact}\label{fac5.3sup}
$\;\;$$\Cl \overline{\phi}, S^* \Cr P^* \rightarrow 
\Cl \overline \phi, S^*P^*  \Cr;$
\end{fact}
so that $P$ preserves $\phi * S$ exactly when 
$P^*$ respects the right lifting $\Cl \overline{\phi}, S^*  \Cr$.\\ 

We now complete the proof of the converse, showing a small
projective weight $\phi$ to be absolute. Supposing $\phi * S$ 
to exist for $S: \B \rightarrow \C$, we are to show that
the right lifting $\Cl \overline{\phi}, S^*  \Cr$ of \ref{fac5.2sup} 
is respected by $P^*$ for every $P: \C \rightarrow \D$.
But this is certainly the case since, $\overline{\phi}$ 
being a left adjoint by Proposition \ref{fladj2}, 
the lifting in question is absolute by Proposition \ref{sladj1}.
\epf
\end{section}

\begin{section}{Cauchy completion and the Morita theorems}\label{section7}
For any category $\A$, the inclusion $J: \A \rightarrow \Q(\A)$  
expresses $Q(\A)$ as the free $Q$-cocomplete category on $\A$,
which by Theorem \ref{Str} is the free cocompletion of $\A$ 
under absolute colimits.
It is determined by the universal property \ref{Pro3.5}, which here,
because every functor preserves {\em absolute} colimits,
becomes: 
\begin{fact}
$\VCat(\Q(\A),\B) \simeq \VCat(\A,\B)$ for any $\B$ with absolute colimits. 
\end{fact}

Proposition \ref{Pro3.15} here takes the following stronger form:
\begin{proposition}
The inclusion $J: \A \rightarrow \Q(\A)$ is an equivalence
if and and only if $\A$ admits all absolute colimits.
\end{proposition}
\pf
The ``only if'' part is trivial. If $\A$ and $\B$ are $\Q$-cocomplete,
we have 
$\QCo(\A,\B) = \VCat(\A,\B) \simeq \QCo(\Q(\A),\B)$,
whence it follows that $J: \A \rightarrow \Q(\A)$ is an equivalence.
\epf

\begin{proposition}\label{pro7.3}
For a small $\B$,
let $\phi \in [\B^{op},\V]$ and $\psi \in [\B,\V]$
be small projective weights related by the equivalence 
\ref{protopro}.
Then for any category $\A$ and any functor 
$F: \B \rightarrow \A$, we have an isomorphism
$\{ \psi, F \} \cong \phi * F$, either side existing 
if the other does. Accordingly, $\A$ admits absolute
limits if and only if it admits absolute colimits.
\end{proposition}
\pf
Let $\{ \psi, F \}$ exist;
as the $\psi$-weighted limit of $F$ in $\A$, it is
also the $\psi$-weighted colimit of $F^{op}$ in $\A^{op}$.
Since $\psi$-weighted colimits are absolute by Theorem \ref{Str},
the canonical 
$ \psi * \A(F-,a) \rightarrow \A( \{ \psi, F \},a)$
is invertible; 
but $\psi * \A(F-,a)$ is isomorphic by \ref{flem} 
to $[\B^{op},\V]( \phi, \A(F-,a))$, exhibiting
$\{ \psi, F \}$ as the colimit $\phi*F$.
\epf

The equivalence \ref{protopro} above was for small categories
$\B$; it admits the following extension to arbitrary categories:
\begin{proposition}\label{pro7.4}
For any category $\A$, we have an equivalence 
${(\Q(\A^{op}))}^{op} \simeq \Q(\A)$.
\end{proposition}
\pf
Let $\B$ admit absolute colimits; so $\B^{op}$ too
admits absolute colimits, by Proposition \ref{pro7.3}.  
Then for any $\A$, we have equivalences
\begin{tabbing}
\spa$\QCo(\Q(\A),\B)$ \=$\simeq$ \=$\VCat(\A,\B)$\\
\>$\cong$ \>${(\VCat(\A^{op},\B^{op}))}^{op}$\\
\>$\simeq$ \>${(\VCat(\Q(\A^{op}),\B^{op}))}^{op}$\\
\>$\cong$ \>$\VCat({(\Q(\A^{op}))}^{op}, \B)$\\
\>$=$ $\QCo( {(\Q(\A^{op})) }^{op}, \B)$.
\end{tabbing}
The desired equivalence ${(\Q(\A^{op}))}^{op} \simeq \Q(\A)$ follows.
\epf

A category $\A$ which admits absolute colimits, and
hence absolute limits, is said to be {\em Cauchy-complete};
and $\Q(\A)$ is called the {\em Cauchy-completion} of $\A$;
this concept was introduced by Lawvere in \cite{Law73}.
For a general class $\Phi$ of weights, the free completion
of $\A$ under {\em $\Phi$-limits} is of course 
${(\Phi(\A^{op}))}^{op}$;
so by Proposition \ref{pro7.4}, $\Q(\A)$ is also the 
completion of $\A$ under {\em absolute limits}.

\begin{proposition}\label{thisprop}
For any class $\Phi$ of weights and any category $\A$, 
the category ${\Phi(\A)}_{\Phi}$ is included in $\Q(\A)$.
If the class $\Phi$ contains $\Q$, we have an equality 
${\Phi(\A)}_{\Phi} = \Q(\A)$.
\end{proposition}
\pf
We begin by proving the first assertion in the case 
of a small $\A$. Let us denote the inclusions again by 
$\xymatrix{ \A  \ar[r]^{Z} & \Phi(\A) \ar[r]^{W} & [\A^{op},\V]}$ 
with $WZ=Y$. For $\phi \in \Phi(\A)_{\Phi}$, the representable
functor $\Phi(\A)(\phi,-):\Phi(\A) \rightarrow \V$ preserves
$\Phi$-colimits; in particular, it preserves the colimit
$\phi * Z \cong \phi$ in $\Phi(\A)$. 
Since $W$ preserves $\Phi$-colimits, 
we have $W(\phi * Z) \cong \phi * WZ = \phi * Y$ 
in $[\A^{op},\V]$.
The composite of 
$[\A^{op},\V](\phi,-) = [\A^{op},\V](W\phi,-)$ 
with $W$ is $[\A^{op},\V](W\phi,W-)$, which is
isomorphic to $\Phi(\A)(\phi,-)$ since $W$ is fully faithful.
Since this composite preserves the colimit $\phi*Z$, it follows
that $[\A^{op},\V](\phi,-)$ preserves the colimit
$W(\phi * Z) \cong \phi*Y$.
Accordingly, $\phi$ is a small projective by \ref{fladj2}.\\

We now prove the first statement for an arbitrary category 
$\A$.  
By Proposition \ref{F3.14}, any $F \in \Phi(\A)$ is of the form
$Lan_{H^{op}}\phi$ for some fully faithful $H: \K \rightarrow \A$
with $\K$ small and some $\phi \in \Phi$. Because $H$ is 
fully faithful, $Lan_{H^{op}}: [\K^{op},\V] \rightarrow \p \A$
is also fully faithful; moreover, as a left adjoint, it 
preserves all colimits. For $\psi \in \Phi(\K)$, its image
$Lan_{H^{op}}\psi = \psi * YH$,
as a $\Phi$-colimit of representables, lies in $\Phi(\A)$;
so that $Lan_{H^{op}}$ restricts to a functor
$L: \Phi(\K) \rightarrow \Phi(\A)$. This functor,
like $Lan_{H^{op}}$, is fully faithful, and it preserves
$\Phi$-colimits, since these are formed in $\Phi(\K)$ as 
in $[\K^{op},\V]$ and in $\Phi(\A)$ as in 
$\p \A$. Since $L$ is fully faithful, we have an
isomorphism $\Phi(\K)(\phi,-) \cong \Phi(\A)(L(\phi),L-) = \Phi(\A)(F,L-)$.
If $F$ belongs to $\Phi(\A)_{\Phi}$ then, $\Phi$-colimits 
are also preserved by $\Phi(\A)(F,-)$, and hence
by $\Phi(\A)(F,L-)$. Thus $\Phi$-colimits are preserved
by $\Phi(\K)(\phi,-)$, so that $\phi$ belongs to $\Q$ 
by the first part of the proof. So $F = \phi * YH$,
as a $\Q$-colimit of representables, lies in $\Q(\A)$.\\

Suppose now that $\Q \subset \Phi$.
Since $\Phi \subset \p$, we have $\p^- \subset \Phi^-$,
or $\Q \subset \Phi^-$. By Proposition \ref{pro5.2}, 
$\Phi(\A)_{\Phi}$ is closed in $\Phi(\A)$ under 
$\Phi^-$-colimits, and hence under $\Q$-colimits.
Since $\Q \subset \Phi$ however, $\Q$-colimits are
preserved by the inclusion $\Phi(\A) \rightarrow \p \A$.
Thus $\Phi(\A)_{\Phi}$ is also closed under $\Q$-colimits 
in $\p \A$; and since it contains the representables,
it contains $\Q(\A)$.
\epf 

For any class $\Phi$ containing $\Q$, for any $\A$ and, for 
any $\Phi$-cocomplete $\B$, we have
$$\PhiCo(\Phi(\A), \B) \simeq \VCAT(\A,\B) \simeq \VCAT(\Q(\A),\B) 
\simeq \PhiCo(\Phi(\Q(\A)), \B),$$
so we have an equivalence,
\begin{fact}\label{morev} 
$\Phi (\A) \simeq \Phi(\Q (\A)).$
\end{fact}
The case $\Phi = \p$ of the following proposition
is the principal classic Morita theorem:
\begin{proposition} Let $\Phi$ be a class of weights
containing $\Q$. Then for any categories $\A$ and $\B$,
we have $\Phi(\A) \simeq \Phi(\B)$ if and 
only if $\Q(\A) \simeq \Q(\B)$.
\end{proposition}
\pf
If $\Q(\A) \simeq \Q(\B)$ we get $\Phi(\A) \simeq \Phi(\B)$
by \ref{morev}. If $\Phi(\A) \simeq \Phi(\B)$ then
$\Q(\A) \simeq \Q(\B)$ by Proposition \ref{thisprop}
\epf
In the circumstances of this proposition,
the categories $\A$ and $\B$ are said to be 
{\em Morita equivalent}.
\end{section}

\begin{section}{$\Phi$-continuous presheaves}\label{section8}
We turn now to the study of categories of the form  
$\PhiCont[ \N^{op}, \V]$, where $\N$ is a small
$\Phi$-cocomplete category and
$\PhiCont[ \N^{op}, \V]$ denotes 
the full subcategory of $[\N^{op},\V]$ 
determined by the $\Phi$-continuous functors 
$\N^{op} \rightarrow \V$.
Since $\PhiCont$ equals $\PhiSCont$ by the definition of 
$\Phi^*$, we may as well suppose that $\Phi$ is
saturated. Since the representables are certainly
$\Phi$-continuous, the Yoneda embedding factorizes 
through the inclusion $J: \PhiCont[\N^{op},\V] \rightarrow [\N^{op},\V]$,
as say $Y =JK$; then since $J$ is fully faithful, it follows from
Yoneda that $J$ is isomorphic (in the notation of Section 
\ref{section4}) to $\tilde{K}$. For the case where $\V = \Set$ and
$\Phi$ is a small set of weights of the form
$\Delta 1: \K^{op} \rightarrow \Set$, some of the results
below appear in \cite{ABLR02}.\\

Since $\Phi$-limits 
commute in $\V$ with all limits and with $\Phi^+$-colimits,
and since such limits and colimits are formed pointwise
in $[\N^{op},\V]$, we have:
\begin{proposition}\label{Pro26}
For any small $\Phi$-cocomplete $\N$, the category 
$\PhiCont[\N^{op},\V]$ is closed in $[\N^{op},\V]$
under all limits and under $\Phi^+$-colimits.
As a consequence, $\Phi^+(\N) \subset \PhiCont[\N^{op},\V]$. 
\end{proposition}
In other words, each 
$\Phi$-flat weight $\N^{op} \rightarrow \V$ is 
$\Phi$-continuous; 
we shall later give conditions for the converse to hold.
That it does not hold in general is shown by the following example,
which was Example 2.3 (vii) of \cite{ABLR02}:
\begin{example}\label{Exa8.2}\end{example}
With $\V = \Set$, let $\Phi$ be the saturated class of weights
for which a $\Phi$-cocomplete category is one with pushouts,
and let $\N^{op}$ be the one-object category given by a non-trivial
group $G$, so that $[\N^{op}, \Set]$ is the category of $G$-sets.
Then $\N^{op}$ has pullbacks, and $\Delta 1: \N^{op} \rightarrow \Set$
preserves pullbacks. Yet $\Delta 1$ is not $\Phi$-flat: for 
\ref{fac2.3} gives $\Delta 1 * - = colim$, which by (3.35) of \cite{Kel82}
sends a presheaf to the set of connected components of its 
set of elements, and thus sends a $G$-set $X$ to the set of 
its orbits. Now the $G$-sets $G \rightarrow 1 \leftarrow G$
have a pullback given by $G \leftarrow G \times G \rightarrow G$,
and this pullback is not preserved by the passage to the sets 
of orbits.\\

Recall from \cite{Kel82} that
many important base-categories $\V$ are {\em locally bounded}, 
and that Theorem $6.11$ of that work gives:
\begin{proposition}\label{Pro32}
Whenever the base category $\V$ is locally bounded,
$\PhiCont[\N^{op},\V]$ is reflective in $[\N^{op},\V]$,
for any class $\Phi$ of weights, and any small 
$\Phi$-cocomplete $\N$. 
\end{proposition}
Sometimes however -- as under certain hypotheses to be introduced 
below -- we can infer the reflectiveness of 
$\PhiCont[\N^{op},\V]$ more easily, without using the general
theorem above, which involves a transfinite induction.
Moreover additional hypotheses may imply special properties of the 
reflexion.\\ 

An important property of $\PhiCont[\N^{op},\V]$ is the following;
this is well known, one generalization of it being 
Theorem 5.56 in \cite{Kel82}.
\begin{lemma}\label{lem29} For any functor $G:\N \rightarrow \B$
where $\N$ is $\Phi$-cocomplete, the corresponding functor
$\tilde{G}: \B \rightarrow [\N^{op},\V]$ takes its 
values in $\PhiCont[\N^{op},\V]$ if and only if $G$ preserves
$\Phi$-colimits.
\end{lemma}
\pf
Consider a $\Phi$-colimit $\phi*T$ in $\N$, where 
$\phi: \LL^{op} \rightarrow \V$ lies in $\Phi$ and 
where $T: \LL \rightarrow \N$.
To say that $G$ preserves this colimit is to say that 
$G( \phi * T)$ (with the appropriate unit) is 
the colimit $\phi * GT$ in $\B$, which is also to say that, 
for each $b \in \B$, the object $\B(\phi* GT,b)$ (with the appropriate
counit) is the limit $\{ \phi, \B(GT-,b) \}$ in $\V$; this, in turn, 
is to say that each $\tilde{G}b: \N^{op} \rightarrow \V$ preserves the 
limit $\{ \phi , T^{op} \}$ in $\N^{op}$. To ask this for each 
$\Phi$-colimit $\phi * T$ in $\N$ is just to ask
$\tilde{G}b$ to lie in $\PhiCont(\N^{op}, \V)$.
\epf

When we take $G: \N \rightarrow \B$ in \ref{lem29}
to be $K: \N \rightarrow \PhiCont[\N^{op},\V]$, it is trivial
that the inclusion
$\tilde{K}: \PhiCont[\N^{op},\V] \rightarrow [\N^{op},\V]$
takes its values in $\PhiCont[\N^{op},\V]$; so the lemma 
gives:
\begin{corollary}\label{Cor30}
For any $\Phi$-cocomplete $\N$, the inclusion 
$K: \N \rightarrow \PhiCont[\N^{op},\V]$
preserves\\$\Phi$-colimits.
\end{corollary}
Another useful lemma is the following:
\begin{lemma}\label{Lem31}
Let $\C$ be a full subcategory of $\A$, 
and write $\B$ for the full subcategory of $\A$ given by
those objects of $\A$ which admit a reflexion 
into $\C$. Let $\C$ and $\A$ admit $\Phi$-colimits.
Then $\B$ is closed in $\A$ under $\Phi$-colimits.
\end{lemma}
\pf
Let the colimit $\phi * T$ exist in $\A$, where 
$\phi: \K^{op} \rightarrow \V$ lies in $\Phi$ and 
where $T: \K \rightarrow \A$ takes its values in $\B$.
Write $R: \B \rightarrow \C$ for the functor sending
each $b \in \B$ to its reflexion $Rb$ in $\C$.
Then for $c \in \C$ we have
\begin{tabbing}
\hspace{1cm}$\A(\phi*T,c)$ \=$\cong$ \=$[\K^{op},\V](\phi, \A(T-,c))$,
by the definition of $\phi*T$;\\
\>$\cong$ \>$[\K^{op},\V](\phi, \C(RT-,c)))$, since $R$ is the
reflexion;\\
\>$\cong$ \>$\C(\phi* RT,c)$, by the definition of $\phi*RT$;
\end{tabbing}
thus $\phi*T$ admits the reflexion $\phi*RT$ in $\C$.
\epf

\begin{remark}\label{Rem33}\end{remark}\noindent
Any full subcategory $\C$ of $[\N^{op},\V]$
is of course cocomplete like $[\N^{op},\V]$ if it is
reflective. It is an old and classical observation
that the converse is also true whenever $\C$ contains
the representables (so that, once again, the inclusion 
$\C \rightarrow [\N^{op},\V]$ is isomorphic to $\tilde{K}$,
where $K: \N \rightarrow \C$ is the inclusion). For to say that 
{\em the object $\phi$ of $[\N^{op},\V]$ admits a reflexion $d$ 
in $\C$} is to say that we have, naturally in $c$, an isomorphism
$$[\N^{op},\V](\phi, \tilde{K}c) \cong \C(d,c),$$
and this is to say that $d$ is the colimit $\phi*K$ in $\C$.\\ 

We shall adopt the following notation: for a full subcategory 
$\B$ of $\A$ and a class $\Phi$ of weights: we write $\Phi \{ \B \}$
for the closure of $\B$ in $\A$ under $\Phi$-colimits.
Of course $\A$ must be understood if this notation
is to suffice: otherwise we should use $\Phi\{\B \mid \A \}$.

\begin{proposition}\label{Pro34}
Still supposing $\N$ to be small and $\Phi$-cocomplete, write 
$\Phi^+ \{ \Phi(\N) \}$ for the closure of $\Phi(\N)$ in 
$[\N^{op},\V]$ under $\Phi^+$-colimits.
Then each object of $[\N^{op},\V]$ that lies in $\Phi^+ \{ \Phi(\N) \}$
has a reflexion in $\PhiCont[\N^{op},\V]$. In fact this reflexion lies
in $\Phi^+(\N)$, and the reflexion of an object of $\Phi(\N)$ lies in $\N$.
\end{proposition}
\pf  Since $\N$ is $\Phi$-cocomplete by hypothesis,
it admits for each $\phi \in \Phi(\N) = \Phi[\N]$
the colimit $\phi * 1_{\N}$ of $1_{\N}: \N \rightarrow \N$.
So since $K: \N \rightarrow \PhiCont[\N^{op},\V]$
preserves $\Phi$-colimits 
by Corollary \ref{Cor30}, the object $K(\phi * 1)$ of 
$\PhiCont[\N^{op},\V]$ is the colimit $\phi * K$; 
and accordingly it is, by Remark \ref{Rem33}, the reflexion in 
$\PhiCont[\N^{op},\V]$ of $\phi \in [\N^{op},\V]$. 
Thus every $\phi \in \Phi(\N)$ has a reflexion in $\PhiCont[\N^{op},\V]$,
which in fact lies in $\N$ (embedded by 
$K: \N \rightarrow \PhiCont[\N^{op},\V]$). Since $\PhiCont[\N^{op},\V]$
admits $\Phi^{+}$-colimits by Proposition \ref{Pro26}, it follows
from Lemma \ref{Lem31} that the objects of $[\N^{op},\V]$ 
admitting a reflexion in $\PhiCont[\N^{op},\V]$ are closed
under $\Phi^+$-colimits; accordingly they include all the objects
of $\Phi^+ \{ \Phi(\N) \}$. Since the reflexion preserves 
$\Phi^+$-colimits, which are (by Proposition \ref{Pro26})
formed in $\PhiCont[\N^{op},\V]$ as they are in $[\N^{op},\V]$, 
the reflexions all lie in the $\Phi^+$-closure $\Phi^+(\N)$ 
of $\N$ in $[\N^{op},\V]$.  
\epf

\begin{theorem}\label{Theo36}
For a small $\Phi$-cocomplete $\N$, the inclusion 
$\Phi^+(\N) \subset \PhiCont[\N^{op},\V]$ of 
Proposition \ref{Pro26} is an equality if and only 
if $\PhiCont[\N^{op},\V] \subset \Phi^+ \{ \Phi(\N) \}$.
In particular $\Phi^+(\N) = \PhiCont[\N^{op},\V]$ whenever
$\Phi^+ \{ \Phi(\N) \}$ is all of $[\N^{op},\V]$.
\end{theorem}
\pf
Since $\N \subset \Phi(\N)$ we have 
$\Phi^+(\N) \subset \Phi^+ \{ \Phi(\N) \}$,
so that certainly $\PhiCont[\N^{op},\V] \subset \Phi^+ \{ \Phi(\N)\}$
if $\Phi^+(\N) = \PhiCont[\N^{op},\V]$.
Now if $\PhiCont[\N^{op},\V] \subset \Phi^+ \{ \Phi(\N)\}$, then by 
Proposition \ref{Pro34} each object of $\PhiCont[\N^{op},\V]$
has for its reflexion in $\PhiCont[\N^{op},\V]$ 
-- namely itself -- an object of $\Phi^+(\N)$.
\epf

We may express the above by saying that, in these
circumstances, the $\Phi$-flat weights coincide with 
the $\Phi$-continuous ones.\\

It is convenient to introduce the following 
definition.
\begin{definition}
A class $\Phi$ of weights is said to be 
{\em locally small} if each $\Phi(\K)$ 
with $\K$ small is also small.
\end{definition}
Since $\Phi^*(\K) = \Phi(\K)$ for any $\K$, 
a class $\Phi$ is locally small if and only if its saturation 
$\Phi^*$ is so.
Moreover, when $\Phi$ is saturated, since we have 
$\Phi(\K) = \Phi[\K]$, to say that $\Phi$ is locally small
is to say that each $\Phi[\K]$ is small. For a general class
$\Phi$, it was observed in Section 3.5 of \cite{Kel82}
that $\Phi$ is locally small when the class $\Phi$ is in fact
a small set. For example, when $\V= \Set$ and $\Phi$ consists
of the three weights giving initial objects, binary coproducts
and coequalizers, $\Phi^*$ is locally small;
here $\Phi(\K)$ is the free finitely-cocomplete category
on $\K$, and $\Phi^*$ is the saturation of the weights
for finite colimits. Similarly when $\V$ is locally finitely
presentable as a closed category, as in \cite{Kel82-2};
what are there called ``the finite indexing types'' form
a small set $\Phi$, so that $\Phi^*$ is locally small; 
here $\Phi(\K)$ is again the free finitely-cocomplete
category on $\K$, and $\Phi^*$ is the saturation of the weights
for finite colimits; compare Examples \ref{exsat} and \ref{exsat1}.
We may note that the class $\Phi$ of example \ref{Exa8.2}
is locally small.\\

For a {\em locally small} saturated class $\Phi$,
the last statement of \ref{Theo36} has a converse.  
First note that for any $K: \A \rightarrow \C$ with $\A$ 
small, the left Kan extension along $K$ of 
$Y: \A \rightarrow [\A^{op},\V]$ is $\tilde{K}$, by 
\ref{fac2.9}. In particular $Lan_Y Y$ is the identity.
However $Y = WZ$, where $Z: \A \rightarrow \Phi(\A)$
and $W: \Phi(\A) \rightarrow [\A^{op},\V]$ are 
the inclusions. Thus 
$Lan_Y Y \cong  Lan_W(Lan_Z Y) \cong Lan_W W$, 
since $Lan_Z Y \cong \tilde{Z} \cong W$. (That is,
in the language of \cite{Kel82}, the functor $W$
is {\em dense}). 
\begin{theorem}\label{teo40}
When the saturated class $\Phi$ is locally small,
the following are equivalent:\\
$(i)$ $\PhiCont[\N^{op},\V] = \Phi^+(\N)$ for any small
$\Phi$-cocomplete $\N$;\\
$(ii)$ For any small $\A$, every presheaf  
$\A^{op} \rightarrow \V$ is a $\Phi^+$-colimit 
of a diagram in $\Phi(\A)$;\\ 
$(iii)$ For any small $\A$, 
$\Phi^{+}\{ \Phi(\A) \} = [\A^{op},\V]$.
\end{theorem}
\pf
$(iii)$ implies $(i)$ by Theorem \ref{Theo36} and
$(ii)$ implies $(iii)$ trivially; so it remains
to prove $(i)$ implies $(ii)$.
For any presheaf $F: \A^{op} \rightarrow \V$,
we have $F \cong Lan_W W(F) \cong [\A^{op},\V](W-,F)* W$
using \ref{fac2.9}.
However the presheaf $[\A^{op},\V](W-,F): \Phi(\A)^{op} \rightarrow \V$
is $\Phi$-continuous as $W$ preserves $\Phi$-colimits
and the representable $[\A^{op},\V](-,F)$ is continuous.
But this presheaf is a weight since $\Phi(\A)$ like $\A$ is
small; so by $(i)$, it belongs to $\Phi^+$.
\epf
\begin{remark} \end{remark} 
A special case of Theorem \ref{teo40} forms part of 
\cite{ABLR02} Theorem 2.4.\\

Notice that $\Phi^+\{ \Phi(\A) \}$ in the theorem above
is different from $\Phi^+(\Phi(\A))$, which is the closure
of $\Phi(\A)$ under $\Phi^+$-colimits in $[\Phi(\A)^{op},\V]$.
Nevertheless since $\Phi \subset \Phi^{+-}$,
$[\A^{op},\V]_{{\Phi}^+}$ is closed in $[\A^{op},\V]$ under 
$\Phi$-colimits by Proposition \ref{pro5.2} and 
since $[\A^{op},\V]_{{\Phi}^+}$ also contains the representables,
it contains $\Phi(\A)$. So Proposition \ref{Pro4.2}
gives:
\begin{observation}\label{Theo39}
Let $\Phi$ be a locally small class satisfying
the equivalent conditions of Theorem \ref{teo40}.
Then for each small $\A$, the category
$[\A^{op},\V]$ is equivalent to 
$\Phi^+(\Phi(\A))$. 
\end{observation}
Note that this result may also be deduced from 
Theorem \ref{Pro3.5} which, because $\V$ is complete,
gives $\PhiCont[\Phi(\A)^{op},\V] \simeq [\A^{op},\V]$
for any small $\A$; for we have 
$\PhiCont[\Phi(\A)^{op},\V] = \Phi^+(\Phi(\A))$ from 
\ref{Theo36}.

\begin{example}\label{Ex41} \end{example}
Let $\V$ be locally finitely presentable as a closed
category, in the sense of \cite{Kel82-2}, and let 
$\Phi$ be the saturation of the weights for finite 
colimits. In this context, the weights in $\Phi^+$
are said to be {\em flat}. For a small $\A$,
the full subcategory $\Phi(\A)$ of $[\A^{op},\V]$
consisting of the finite colimits of the representables 
is also, as mentioned previously in \ref{exsat},
the subcategory of 
$[\A^{op},\V]$ given by the finitely presentable
objects; moreover, by \cite{Kel82-2} Theorem 7.2 again, 
every object of $[\A^{op},\V]$
is a filtered colimit of a diagram in $\Phi(\A)$. However
(conical) filtered colimits in $\V$ commute,
by \cite{Kel82-2} Proposition 4.9, with finite limits; 
so the weight for a conical filtered colimit is flat -- that is,
belongs to
$\Phi^+$. Thus $\Phi^+ \{ \Phi(\A) \}$ is all of $[\A^{op},\V]$
for any small $\A$, so that Theorem \ref{teo40} applies
for this $\Phi$.

\begin{example}\label{Ex42}\end{example}
Everything in Example \ref{Ex41} continues to hold
when we take for $\Phi$ not the weights for finite 
colimits but those for $\alpha$-colimits, where 
$\alpha$ is a regular cardinal; see \cite{Kel82-2},
Section 7.4.

\begin{example}\label{Ex45}\end{example}
Let $\Phi$ be the saturation of the class $\Psi$
of Example \ref{exsat2}, so that a $\Phi$-cocomplete
category is one with an initial object.
Here $\Phi(\A)$ consists of the representables
along with the initial object $0$ of $[\A^{op},\Set]$. 
Now any presheaf $F :\A^{op} \rightarrow \Set$
is the conical colimit 
of the canonical $Y / F \rightarrow [\A^{op},\Set]$, where 
$Y / F$ is the comma category of 
$Y: \A \rightarrow [\A^{op},\Set]$
and $F: 1 \rightarrow [\A^{op}, \Set]$; 
so it is also the conical colimit of the canonical
$W/F \rightarrow [\A^{op},\Set]$, where 
$W: \Phi(\A) \rightarrow [\A^{op},\Set]$ is the inclusion; 
for $W/F$ differs from $Y/F$ only by the addition
of an initial object, namely the unique map $0 \rightarrow F$.
Since $W/F$ is accordingly connected,
$F$ lies in $\Phi^+ \{ \Phi(\A) \}$ by Example \ref{exsat2}. 
So again Theorem \ref{teo40} applies.

\begin{remark}\label{Ex43}\end{remark}
We get a trivial case where Theorem \ref{teo40} applies
if we take $\Phi$ to be $\Q$. Since $\Q = \p^- = 0^{+\;-}$,
we have $\Q^+ = 0^{+\;-\;+} = 0^+ = \p$.
By \cite{Joh89}, the class $\Q$ is locally small if 
$\V_0$ is locally presentable.
\end{section}

\bibliographystyle{alpha}

\begin{thebibliography}{}
\bibitem[ABLR02]{ABLR02}
{\sc J.Ad\`{a}mek, F.Borceux, S.Lack, and J.Rosick\'{y}},
{\it A classification of accessible categories},
J Pure Appl. Algebra 175 (2002) 7-30.
\bibitem[AR94]{AR94}
{\sc J.Ad\`{a}mek and  J.Rosick\'{y}},
{\it Locally presentable and accessible categories},
Cambridge University Press, (1994).
\bibitem[AK88]{AK88}
{\sc M.H.Albert and G.M.Kelly},
{\it The closure of a class of colimits},
J Pure Appl. Algebra 51 (1988) 1-17.
\bibitem[Bet85]{Bet85} {\sc R.Betti},
{\it Cocompleteness over coverings},
J. Austral. Math. Soc. (Series A) 39 (1985) 169-177.
\bibitem[BQ96]{BQ96}
{\sc F.Borceux and C.Quintero},
{\it Enriched accessible categories},
Bull Austral Math Soc 54 (1996) 489-501.
\bibitem[BQR98]{BQR98}
{\sc F.Borceux, C.Quintero, and J.Rosick\'{y}},
{\it A theory of enriched sketches},
Theory and Applications of Categories 4 (1998) 47-72.
\bibitem[Joh89]{Joh89}
{\sc S.R.Johnson},
{\it Small Cauchy completions},
J. Pure Appl. Algebra 62 (1989) 35-45.
\bibitem[KeLa80]{KeLa80}
{\sc G.M.Kelly and M.Laplaza},
{\it Compact closed categories},
J. Pure Appl. Algebra 19 (1980) 193-213.    
\bibitem[Kel82]{Kel82}
{\sc G.M.Kelly},
{\it Basic concepts of enriched category theory},
London Mathematical Society Lecture Note Series 64,
Cambridge University Press (1982).
\bibitem[Kel82-2]{Kel82-2}
{\sc G.M.Kelly},
{\it Structures defined by finite limits in the enriched context, I},
Cahiers de Top. et G\'eom. Diff. 23 (1982) 3-42.
\bibitem[Law73]{Law73}
{\sc F.W.Lawvere},
{\it Metric spaces, generalized logic, and closed categories},
Rend. Sem. Mat. Fis. di Milano 43 (1973) 135-166,
Reprint in Theory and Applications of Categories, No. 1, 
(2002) 1-37.
\bibitem[Lin74]{Lin74}
{\sc H.Lindner}
{\it Morita equivalences of enriched categories},
Cahiers de Top. et G\'eom. Diff. 15 (1974) 377-397.
\bibitem[MP89]{MP89}
{\sc M.Makkai and R.Par\'e},
{\it Accessible categories: the foundations of categorical 
model theory},
Contemporary Math. 104, American Math. Soc., Providence (1989).
\bibitem[Str83]{Str83}
{\sc R.Street},
{\it Absolute colimits in enriched categories},
Cahiers de Top. et G\'eom. Diff. 24 (1983) 377-379.
\bibitem[StWa78]{StWa78}
{\sc R.Street and R.Walters},
{\it Yoneda structures on 2-categories},
J. Algebra 50 (1978) 350-379.
\end{thebibliography}

\end{document}